\documentclass[11pt]{article}

\usepackage{amsmath}  
\usepackage{amsthm} 
\usepackage{amssymb} 
\usepackage{hyperref} 
\usepackage{enumitem} 

\usepackage{color} 

\newtheorem{thm}{Theorem}
\newtheorem{inspr}[thm]{}

\newenvironment{definitie}{\begin{itemize}\item[ ]\hspace{-26pt}\bf Definition \rm }{\end{itemize}}
\newenvironment{notatie}{\begin{itemize}\item[ ]\hspace{-26pt}\bf Notation \rm }{\end{itemize}}
\newenvironment{voorbeeld}{\begin{itemize}\item[ ]\hspace{-26pt}\bf Example \rm }{\end{itemize}}
\newenvironment{stelling}{\begin{itemize}\item[ ]\hspace{-26pt}\bf Theorem \rm }{\end{itemize}}
\newenvironment{propositie}{\begin{itemize}\item[ ]\hspace{-26pt}\bf Proposition \rm }{\end{itemize}}
\newenvironment{lemma}{\begin{itemize}\item[ ]\hspace{-26pt}\bf Lemma \rm }{\end{itemize}}
\newenvironment{opmerking}{\begin{itemize}\item[ ]\hspace{-26pt}\bf Remark \rm }{\end{itemize}}
\newenvironment{voorwaarde}{\begin{itemize}\item[ ]\hspace{-26pt}\bf Condition \rm }{\end{itemize}}
\newenvironment{probleem}{\begin{itemize}\item[ ]\hspace{-26pt}\bf Problem \rm }{\end{itemize}}
\newenvironment{gevolg}{\begin{itemize}\item[ ]\hspace{-26pt}\bf Corollary \rm }{\end{itemize}}
\newenvironment{niets}{\begin{itemize}\item[ ]\hspace{-26pt}\bf   \rm }{\end{itemize}}

\newcommand{\defin}{\begin{inspr}\begin{definitie}}  
\newcommand{\edefin}{\end{definitie}\end{inspr}}

\newcommand{\notat}{\begin{inspr}\begin{notatie}}  
\newcommand{\enotat}{\end{notatie}\end{inspr}}

\newcommand{\voorb}{\begin{inspr}\begin{voorbeeld}}  
\newcommand{\evoorb}{\end{voorbeeld}\end{inspr}}

\newcommand{\stel}{\begin{inspr}\begin{stelling}}
\newcommand{\estel}{\end{stelling}\end{inspr}}

\newcommand{\prop}{\begin{inspr}\begin{propositie}}
\newcommand{\eprop}{\end{propositie}\end{inspr}}

\newcommand{\lem}{\begin{inspr}\begin{lemma}}
\newcommand{\elem}{\end{lemma}\end{inspr}}

\newcommand{\opm}{\begin{inspr}\begin{opmerking}}
\newcommand{\eopm}{\end{opmerking}\end{inspr}}

\newcommand{\voorw}{\begin{inspr}\begin{voorwaarde}}
\newcommand{\evoorw}{\end{voorwaarde}\end{inspr}}

\newcommand{\probl}{\begin{inspr}\begin{probleem}}
\newcommand{\eprobl}{\end{probleem}\end{inspr}}

\newcommand{\gev}{\begin{inspr}\begin{gevolg}}
\newcommand{\egev}{\end{gevolg}\end{inspr}}

\newcommand{\nul}{\begin{inspr}\begin{niets}}
\newcommand{\enul}{\end{niets}\end{inspr}}

\newcommand{\bew}{\vspace{-0.3cm}\begin{itemize}\item[ ] \bf Proof\rm: }
\newcommand{\ebew}{\hfill $\qed$ \end{itemize}}

\newcommand{\ssnl}{\vskip 3pt} 
\newcommand{\snl}{\vskip 7pt} 
\newcommand{\nl}{\vskip 12pt} 

\newcommand{\ot}{\otimes}

\newcommand{\inv}{^{-1}}

\newcommand{\tl}{\triangleleft}

\newcommand{\tussenen}{\qquad\quad\text{and}\qquad\quad}

\newcommand{\rood}{\color{red}}
\newcommand{\blauw}{\color{blue}}

\numberwithin{thm}{section}   
\numberwithin{equation}{section} 

\newcommand{\keepcomment}[1]{}
\newcommand{\oldcomment}[1]{}

\begin{document}

\centerline{\bf \Large Polynomial functions for locally compact group actions}
\vspace{13pt}
\centerline{\it Magnus B. Landstad \rm $^{(1)}$ and \it Alfons Van Daele \rm $^{(2)}$}
\bigskip\bigskip

{\bf Abstract} 
\nl
Consider a locally compact group $G$ and a locally compact space $X$. A \emph{local right action} of $G$ on $X$ is a continuous map $(x,p)\mapsto x\cdot p$ from an open subset $\Gamma$ of the Cartesian product $X\times G$ to $X$ satisfying certain obvious properties. 
A global right action of $G$ on $X$ gives rise to a global left action of $G$ on the space $C_c(X)$ of continuous complex functions with compact support in $X$ by the formula $p\,\cdot f:x\mapsto f(x\cdot p)$. In the case of a local action, one still can define $p\,\cdot f$ in $C_c(X)$ by this formula for $f\in C_c(X)$ and $p$ in a neighborhood $V_f$ of the identity in $G$. This yields a local left action of $G$ on $C_c(X)$.
\ssnl
Given a local right action of $G$ on $X$, a function $f\in C_c(X)$ is called \emph{polynomial} if there is a neighborhood $V$ of the identity, contained in $V_f$, and a finite-dimensional subspace $F$ of $C_c(X)$ containing all the functions $v\cdot f$ for $v\in V$. In this paper we study such polynomial functions.
\ssnl
If $G$ acts on itself by multiplication, we are also interested in the local actions obtained by restricting it to an open subset of $G$. This is the typical situation that is encountered in our paper on bicrossproducts of groups with a compact open subgroup, \cite{La-VD3}. In fact, the need for a better understanding of polynomial functions for that case has led us to develop the theory in general here. 

\nl
Date: {\it 15 September 2023}
\nl
\vskip 4 cm
\hrule
\medskip
\begin{itemize}
\item[$^{(1)}$] Department of Mathematical Sciences, Norwegian University of Science  and  Technology, NO-7491 Trondheim, Norway. E-mail: \texttt{magnus.landstad@ntnu.no}
\item[$^{(2)}$] Department of Mathematics, University of Leuven, Celestijnenlaan 200B,
B-3001 Heverlee, Belgium. E-mail: \texttt{alfons.vandaele@kuleuven.be}
\end{itemize}
\newpage

\setcounter{section}{-1}  

\section{\hspace{-17pt}. Introduction} \label{s:intro}  

Let $G$ be a locally compact group and $X$ a locally compact space. To define a local right action of $G$ on $X$ we first need an open subset $\Gamma$ of the Cartesian product $X\times G$ and further a continuous map $(x, p)\mapsto x\cdot p$ from this open subset $\Gamma$ to $X$. We require some natural conditions. In the first place we want $(x,e)\in \Gamma$ and $x\cdot e=x$ for all $x\in X$. Here $e$ denotes the identity in $G$. Further, given $(x,p)\in \Gamma$ and any other element $q$ in $G$, we want that $(x,pq)\in \Gamma$ if and only if $(x\cdot p,q)\in \Gamma$ and then $x\cdot pq=(x\cdot p)\cdot q$.
\ssnl
A local right action with $\Gamma=X\times G$ is a global right action. On the other hand, given a global right action of $G$ on $X$ and an open subset $Y$ of $X$ we have a natural local right action of $G$ on $Y$ if we let $\Gamma$ be the set of pairs $(y,p)$ where $y\in Y$ and also $y\cdot p\in Y$. 
It is obvious to show that this yields a local right action on $Y$. 
\ssnl
One of the questions one may ask is whether any local  action is the restriction of a global  action as above. The answer is probably no, but many interesting local actions  are of that form. Therefore, properties of global actions play an important role for studying local actions. 
\ssnl
Given a global right action of $G$ on $X$, we have an obvious  global left action of $G$ on the space $C_c(X)$ of continuous complex functions with compact support in $X$. It is given by $p\cdot f:x\mapsto f(x\cdot p)$ for $x\in X$ and $p\in G$. The continuity of the action of $G$ on $X$ implies that $p\cdot f$ is again a function in $C_c(X)$. The situation is completely different in the case of a local action. Now $(x,p)\mapsto f(x\cdot p)$ is defined and continuous on $\Gamma$. One can define $p\cdot f$ on all of  $X$ by $(p\cdot f)(x)=f(x\cdot p)$ when $(x,p)\in \Gamma$ and $(p\cdot f)(x)=0$ otherwise. In general, this function will no longer be continuous. Fortunately, given $f\in C_c(X)$ there exists a neighborhood $V_f$ of the identity so that the  function $v\cdot f$ for $v\in V_f$, as defined above, is continuous and again belongs to $C_c(X)$. In this way, we obtain a local left action of $G$ on $C_c(X)$.
\ssnl
Now polynomial functions can be defined. We call $f\in C_c(X)$ polynomial if there is a neighborhood $V$  of $e$ contained in $V_f$, and a finite-dimensional subspace $F$ of $C_c(X)$ containing all the functions $v\cdot f$ for $v\in V$. 
\ssnl
Polynomial functions for the action of $G$ on itself have been studied first \keepcomment{true?}in \cite{La-VD1}. They appear naturally when the existence of a multiplier Hopf algebra in the function algebra $C_c(G)$ is investigated. This is shown, among other things, in \cite{La-VD1}. In \cite{La-VD3} the same question is posed for the bicrossproduct construction associated to an admissible pair of closed subgroups of a group $G$. It takes us quickly to the situation where $G$ has a compact open subgroup and to algebras of polynomial functions. This more general case is far more complicated and a deep knowledge of polynomial functions for local actions is necessary.
\ssnl
This is the motivation for writing this paper, devoted to the study of such polynomial functions for local actions. The reader finds interesting applications of the theory developed in this paper for the case of an admissible pair of closed subgroups of a locally compact group in \cite{La-VD3}. We have included references here to the particular results in \cite{La-VD3}.
\nl
\bf Content of the paper \rm
\nl
In \emph{Section} \ref{s:loc-act} we define and study local actions of a locally compact group $G$ on a locally compact space $X$. We develop the theory mostly for local right  actions, but also left actions and commuting left and right actions are considered. There are two special properties of actions that play an important role. First a \emph{proper action} is an action where the map $(x,p)\mapsto (x,x\cdot p)$ from $\Gamma$ to $X\times X$ is a proper map. The other one is where, for each $x$ in $X$, we have a neighborhood $V_x$ of the identity in $G$ with the property that $(x,v)\in \Gamma$ for all $v\in V_x$ and such that $v\mapsto x\cdot v$ is a homeomorphism from $V$ to its image in $X$. We say that the action is \emph{locally homeomorphic}.
\ssnl
In this section we also consider different methods to obtain new actions from a given one.
\snl
In \emph{Section} \ref{s:ind-fun} we define the associated local left action of the group on the space $C_c(X)$. We show that for each $f\in C_c(X)$ there is a neighborhood $V_f$ of the identity in $G$ so that $v\cdot f$ is well-defined in $C_c(G)$ for all $v\in V_f$. We do not give a formal definition but derive some properties that would be expected for a local action of the group on $C_c(X)$.
\snl
This short section allows us to define polynomial functions in \emph{Section} \ref{s:pol} where we derive its first properties. In particular, given a polynomial function $f$ with compact support, there is a compact set $C$ in $X$, containing the support of $f$, and an open and closed subgroup $G_0$ of $G$ so that $(x,p)\in\Gamma$ and $x\cdot p\in C$ for all $x\in C$ and $p\in G_0$. In the case of a proper action, this group $G_0$ will be a compact open subgroup of $G$. 
\snl
This result takes us to \emph{Section} \ref{s:cpt-open}. Here we assume that the group $G$ has a compact open subgroup $G_0$. We consider compact subsets $C$ of $X$ with the property as above, namely that $C\cdot G_0$ is well-defined and equal to $C$. The space $\mathcal F_C$, defined as the subspace of $C_c(X)$ of functions with compact support in $C$, is invariant under the action of $G_0$. We get a representation of the compact group $G_0$. With the help of the representation theory of compact groups, we can now construct polynomial functions.
\snl
\emph{Section} \ref{s:pol-lcg} is devoted to the particular case of a group acting on itself by right multiplication. Polynomial functions can only exist when $G$ has a compact open subgroup. They are characterized in $C_c(G)$ and a general method is used to construct polynomial functions. This is mostly contained already \cite{La-VD1}, but we also give some more results. Among other things, we show that for both the function algebra and the convolution algebra, $C_c(G)$ contains the polynomial functions as a $^*$-subalgebra with local units. This gives some information not present in  \cite{La-VD1}.
\snl
In \emph{Section} \ref{s:misc} we treat various properties of polynomial functions in different situations. We show e.g.\ that under certain circumstances, a polynomial function for the action of a subgroup will automatically be also polynomial for the original action of the bigger group. We include a couple of interesting applications of that result.
\snl
Finally in Section \ref{s:concl} we look at possible further research and draw some general conclusions.
\nl 
\bf Notations and conventions \rm
\nl
When $X$ is a locally compact space, we use $C_c(X)$ for the vector space of continuous complex functions on $X$ with compact support. If $X$ is compact we simply use $C(X)$ for the  space of continuous complex functions on $X$. By the support of a function, we always mean the \emph{closure} of the set of points where the function is non-zero. 
\ssnl 
We use $e$ for the identity in a group. By a neighborhood of the identity in a locally compact group, we always mean an \emph{open} neighborhood. 
When $G$ is a group, we use $G^{\text{op}}$ for group obtained from $G$ by taking the opposite multiplication. 
\oldcomment{\ssnl We may need more things here. To be completed ?}{}
\newpage
\bf Basic references \rm
\nl
For the theory of locally compact groups, in particular the representation theory of compact groups, our main references are  \cite{Fo, He-Ro}. Further we refer to \cite{La-VD1} where polynomial functions on a locally compact group have been studied. Finally, \cite{La-VD2} and \cite{La-VD3} are references for some of the examples we use in this note to illustrate the main results.
\oldcomment{\ssnl To be completed ?}
\nl
\bf Acknowledgments \rm
\nl 
We want to thank NTNU (Norway) and KU Leuven (Belgium) for the opportunities we
get to continue our research after our official retirements.

\section{\hspace{-17pt}. Local actions} \label{s:loc-act}  

Assume that $G$ is a locally compact group, acting \emph{globally} on a locally compact space $X$. It is interesting to consider the behavior of such an action on an open subset $Y$ of $X$.  In our papers on bicrossproducts \cite{La-VD2,La-VD3}, we encounter examples of this type. As explained already in the introduction, these examples  have led us to the notion of a \emph{local action} and the study of polynomial functions in this setting. 
\nl
\bf Local actions of a locally compact group on a locally compact space \rm
\nl
Here is the definition we will work with. 

\defin\label{defin:1.1}
Let $G$ be a locally compact group and $X$ a locally compact space. A \emph{local right action} of $G$ on $X$ is a continuous map $(x,p)\mapsto x\cdot p$ from an open set $\Gamma$ of the Cartesian product $X\times G$ to $X$ satisfying the following two conditions:
\vspace* {-5pt}
\begin{itemize}[noitemsep]
\item [i)] The pair $(x,e)$ belongs to $\Gamma$  and $x\cdot e=x$ for all $x\in X$. 
\item[ii)] Given a pair $(x,p)\in \Gamma$ and any other element $q\in G$, then $(x,pq)\in \Gamma$ if and
only if $(x\cdot p, q)\in \Gamma$. In that case $x\cdot (pq)=(x\cdot p)\cdot q$.
\end{itemize}
\edefin

It is clear that a local action is a global action if and only $\Gamma=X\times G$. In that case, the first statements in item i) and ii) of the definition above are obvious and it is just required that $x\cdot e=x$ and $x\cdot pq=(x\cdot p)\cdot q$ for all $x\in X$ and $p,q\in G$.
\ssnl
In what follows, when we say that $x\cdot p$ is defined, we mean that $(x,p)\in\Gamma$. Sometimes formulations will be simpler with this terminology. 
\ssnl
In a similar way we define \emph{local left actions}. In this case we consider continuous maps $(p,x)\mapsto p\cdot x$ from an open subset of $G\times X$ to $X$. 
\ssnl
One can think of local actions in terms of a groupoid as we show in the following proposition.

\prop\label{prop:1.2}
Let $(x,p) \mapsto x\cdot p$ be a local right action of a locally compact group $G$ on a local compact space $X$, defined on the open subset $\Gamma$ of $X\times G$ as in Definition \ref{defin:1.1}. Then $\Gamma$ is a topological groupoid for  the product 
$(x_1,p_1)(x_2,p_2)=(x_1,p_1p_2)$ if 
$x_1\cdot p_1=x_2$. 
The set of units can be identitied with $X$ by the map $x\mapsto (x,e)$. Then the target map is $(x,p)\mapsto x$ while the source map is $(x,p)\mapsto x\cdot p$.
\eprop
\bew
i) The product is well-defined because of condition ii) in Definition \ref{defin:1.1}. To prove the coassociativity of the product, we essentially need that, given  elements $(x_1,p_1)$, $(x_2,p_2)$ and $(x_3,p_3)$ in $\Gamma$, if $x_2=x_1\cdot p_1$ then 
\begin{equation*}
x_2\cdot p_2=x_3 \quad\text{if and only if}\quad x_1\cdot (p_1p_2)=x_3.
\end{equation*}
The coasociativity of the product then is again a consequence of this condition in the definition and of the  coassociativity of the product in the group.
\ssnl
ii) For all $(x,p)$ in $\Gamma$  we have
\begin{equation*}
(x,e)(x,p)=(x,p) \tussenen (x,p)(x\cdot p,e)=(x,p).
\end{equation*}
\ssnl
iii) Finally we have for all $(x,p)\in \Gamma$,
\begin{equation*}
(x,p)(x\cdot p\inv)=(x,e)
\tussenen
(x\cdot p,p\inv)(x,p)=(x\cdot p,e).
\end{equation*}
\ssnl
iv) The range and source maps, the inverse and the product are continuous maps. 

\ebew

Next we have a notion of \emph{commuting actions}. The definition is a bit more involved than it would be for global actions. First consider a combination of a left and a right action.

\defin\label{defin:1.2}
Let $X$ be a locally compact space and $H$ and $K$  locally compact groups. Suppose that we have a local left action $(h,x)\mapsto h\cdot x$, defined from the open set $\Gamma_H$ of $H\times X$ to $X$ and a local right action $(x,k)\mapsto x\cdot k$ from the open set $\Gamma_K$ of $X\times K$ to $X$. Given $x\in X$, $h\in H$ and $k\in K$ satisfying $(h,x)\in\Gamma_H$ and $(x,k)\in\Gamma_K$, 
  we assume that 
  \begin{equation*}
(h\cdot x, k)\in \Gamma_K \text{ if and only if }
(h,x\cdot k)\in\Gamma_H.
\end{equation*}
If these conditions are satisfied, we require that $(h\cdot x)\cdot k=h\cdot (x\cdot k)$. We say that the  local left action of $H$ and the  local right action of $K$ \emph{commute}.
\edefin

In the case of global actions, we simply have a left action of $H$ and a right action of $K$ satisfying $(h\cdot x)\cdot k=h\cdot (x\cdot k)$ for all elements $x,h,k$. For local actions, we need to be more precise. We can say that $(h\cdot x)\cdot k$ is defined if and only if $h\cdot (x\cdot k)$ is defined and that then these elements are equal.
\ssnl
In a similar way we can define commuting  local left actions and commuting  local right actions. 

\prop\label{prop:1.3}
Suppose that we have commuting local right actions of the groups $H$ and $K$ on $X$. Then we have a local right action of the Cartesian product $H\times K$ given by $(x,(h,k))\mapsto (x\cdot h)\cdot k$ when $(x,h)\in\Gamma_H$ and $(x\cdot h,k)\in\Gamma_K$. 
\eprop
\bew
Define $\Gamma$ as the set of elements $(x,(h,k))$ in $X\times(H\times K)$ satisfying $(x,h)\in \Gamma_H$, $(x,k)\in \Gamma_K$ and $(x\cdot h,k)\in \Gamma_K$. Then, by assumption, also $(x\cdot k, h)\in \Gamma_H$ and $(x\cdot h)\cdot k=(x\cdot k)\cdot k$. Further we define $x\cdot (h,k)=(x\cdot h)\cdot k$ for $(x,(h,k))\in\Gamma$. So we also have $x\cdot (h,k)=(x\cdot k)\cdot h$.
\ssnl
This property is essential to show that the action of $H\times K$ satisfies condition ii) of Definition \ref{defin:1.1}. But the argument is straightforward. 
\ebew

Conversely, assume that we are given a right local action of the group $H\times K$ on $X$, defined on an open subset $\Gamma$ of $X\times (H\times K)$. It seems reasonable to assume that $(x,(h,k))\in\Gamma$ implies that both $(x,(h,e))$ and $(x,(e,k))$ belong to $\Gamma$. By condition ii) in Definition \ref{defin:1.1}, we will have that 
$(x,(h,k))\in\Gamma$ if and only if $(x,(h,e))$ and $(x,(e,k))$ belong to $\Gamma$ and 
\begin{equation*}
(x\cdot (h,k))=(x\cdot h)\cdot k=(x\cdot k)\cdot h.
\end{equation*}
So we have to define  $\Gamma_H$ as the sef of elements $(x,h)$ satisfying $(x,(h,e))\in \Gamma$ and $\Gamma_K$ as the set of elements $(x,k)$ satisfying $(x,(e,k))\in \Gamma$.

\keepcomment{\rood \ssnl We have to check this item once more.}{}
\ssnl
If we have a  local left action of $H$, we have an obvious local right action of $H^\text{op}$. If the left action of $H$ commutes with the right action of $K$, then the local right actions of $H^\text{op}$ and $K$ will commute.  The above proposition can be applied and yields a local right action of $H^\text{op}\times K$ on $X$.
\ssnl
We will encounter examples of this later in the paper. We have e.g.\ Proposition \ref{prop:5.12} where we consider the actions of $G$ on itself, given by left multiplication and right multiplication. Further we have Proposition \ref{prop:6.17a}, see also Remark \ref{opm:6.17}.
\ssnl
The following map will play a role further. It is intimately related with the condition in item ii) of the Definition \ref{defin:1.1} above. 

\prop\label{prop:1.4}
Assume that we have a right local action as in Definition \ref{defin:1.1}. Define $\gamma$ from $\Gamma$ to $X\times G$ by $\gamma(x,p)=(x\cdot p,p\inv)$. It is a continuous  bijective map from $\Gamma$ to itself. The inverse is again $\gamma$ and hence it is a homeomorphism. 
\eprop

\bew
i) The map is well-defined and continuous because $(x,p)\mapsto x\cdot p$ is continuous from $\Gamma$ to $X$ and $p\mapsto p\inv$ is continuous from $G$ to itself. 
\ssnl
ii) Take any pair $(x,p)\in\Gamma$. By assumption we have $(x,e)\in \Gamma$ and $x\cdot e=x$. As $pp\inv=e$, from item ii) of the definition we get that $(x\cdot p,p\inv)\in \Gamma$ and that $(x\cdot p)\cdot p\inv =x\cdot e=x$. Hence $\gamma(x,p)\in \Gamma$ and
\begin{equation*}
\gamma(\gamma(x,p))=\gamma(x\cdot p,p\inv)=((x\cdot p)\cdot p\inv, p)=(x,p).
\end{equation*} 
Then the result follows.
\ebew

Let $X_p$ be the set of elements $x\in X$ satisfying $(x,p)\in\Gamma$ for a given element $p\in G$. Then we see from the above result that the map $x\mapsto x\cdot p$ is a homeomorphism from $X_p$ to $X_{p\inv}.$
\oldcomment{\rood Add here a statement about $X_p$}{}

\nl
\bf Proper actions and locally homeomorphic actions \rm
\nl
There are two properties of an action that will play an important role for our study of polynomial functions.

\defin\label{defin:1.5}
A local right action is called \emph{proper} if the map $(x,p)\mapsto (x,x\cdot p)$ from $\Gamma$ to $X\times X$ is a proper map, i.e.\ if the inverse image of a compact subset of $X\times X$ is a compact subset of $\Gamma$.
\edefin

This is exactly the same definition as for global actions. Various notions of properness of actions have been studied in \cite{K-L-Q}. 
\ssnl
The right global action of $G$ on itself, given by right multiplication, is a proper action because the map $(p,q)\mapsto (p,pq)$ has a continuous inverse $(p,q)\mapsto (p,p\inv q)$.
\nl
The second property of a local action we will  consider is the following.

\defin\label{defin:1.6}
A local right action is called \emph{locally homeomorphic} if for every $x\in X$ there exists a neighborhood $V_x$ of the identity in $G$ so that 
\vspace* {-5pt}
\begin{itemize}[noitemsep]
\item[i)] $(x,p)\in\Gamma$ for all $p\in V_x$,
\item[ii)] the set $x\cdot V_x$ of points $x\cdot p$ with $p\in V_x$ is open in $X$,
\item[iii)] the map $p\mapsto x\cdot p$ is a homeomorphism from $V_x$ to $x\cdot V_x$. 
\end{itemize}
We  use $p_x$ for the inverse of the map in item iii). So $x\cdot p_x(y)=y$ for $y\in x\cdot V_x$.
\edefin

Note that if we replace $V_x$ by a smaller neighborhood, we will still have the
same property, precisely because of iii).
\ssnl
Observe also that the right global action of $G$ on itself, given by right multiplication, satisfies this property. In this case, $V_x$ is all of $G$ and $x\cdot V_x=G$ for each $x\in G$. The map $p_x$ is defined on all of $G$ and $p_x(y)=x\inv y$. 
\ssnl
Consider the groupoid structure on $\Gamma$ as in Proposition \ref{prop:1.2}. Recall that a topological groupoid is called \emph{\'etale} (see \cite{Re,Si} if the target map is a local homeomorphism. For the groupoid associated to a local right action, the target map is given by $(x,p)\mapsto x$. This can only be a local homeomorphism if the group $G$ is discrete.  We can conclude from this that the property defined above, although similar to the condition for an \'etale groupoid, is not related in any natural way.
\ssnl
The notion of a locally homeomorphic action is, as far as we know, not a common property, even for global actions.
\nl
\bf Restriction of actions \rm
\nl
Now we show that our \emph{motivating example} satisfies the conditions of a local action. We obtain it as a special case of the following more general construction.

\prop\label{prop:1.7}
Assume that we are given a local right action of $G$ on $X$ as in Definition \ref{defin:1.1}. Let $Y$ be an open subset of $X$ and write
\begin{equation*}
\Gamma_Y=\left\{(y,p)\in\Gamma\mid y\in Y \text{ and } y\cdot p\in Y\right\}.
\end{equation*}
Then we have a  local right action of $G$ on $Y$ given by $(y,p)\mapsto y\cdot p$ where $(y,p)\in \Gamma_Y$.
\eprop
\bew
i)  Because the original action is a continuous map from $\Gamma$ to $X$ and because $Y$ is an open subset of $X$, we have that $\Gamma_Y$ is an open subset of $Y\times G$. For the action of $G$ on $Y$ we just take the restriction of the original action to $\Gamma_Y$. It will be a continuous map.
\ssnl
ii) Clearly for all $y\in Y$ we have $(y,e)\in \Gamma_Y$ and $y\cdot e=y$. 
\ssnl
iii) Now assume that $(y,p)\in \Gamma_Y$ and let $q\in G$. For the original action we know that $(y\cdot p,q)\in \Gamma$ if and only if $(y,pq)\in \Gamma$ and that then $(y\cdot p)\cdot q=y\cdot pq$. Now we have that $y\cdot p\in Y$ and $(y\cdot p)\cdot q\in Y$ if and only if $y\cdot pq\in Y$. Hence we get the required compatibility property also for the action of $G$ on $Y$.
\ebew

We will call the action of $G$ on $Y$ obtained like this simply \emph{the restriction to $Y$} of the original action of $G$ on $X$. Observe that we use the same notation for the action and for its restriction to $Y$.
\ssnl
We will  give two examples to illustrate this procedure. But first, we prove a few properties.
\ssnl
As a special but important case of Propostion \ref{prop:1.7}, we have the following result.

\prop\label{prop:1.8}
Assume that $G$ is a locally compact group, acting globally on a locally compact space $X$. Consider a non-empty open subset $Y$ of $X$. Define $\Gamma$ as the set of pairs $(y,p)$ where $y\in Y$ and $p\in G$ satisfying $y\cdot p\in Y$.  The map $(y,p)\mapsto y\cdot p$ from $\Gamma$ to $Y$ is a local action in the sense of Definition \ref{defin:1.1}.
\eprop

Remark that if the original action is proper, the same is true for the associated local action on the subspace. 
For locally homeomorphic actions we have a similar result.

\prop\label{prop:1.9}
Suppose that we have a local right action of $G$ on $X$  that it is locally homeomorphic. Let $Y$ be an open subset of $X$. Then the  local right action of $G$ on $Y$, obtained by restricting the global action is still locally homeomorphic.
\eprop

\bew
For $x\in X$, denote by $V_x$ a neighborhood of the identity in $G$ as in the Definition \ref{defin:1.6} above. When $y\in Y$ we put
\begin{equation*}
V'_y=\{p\in G\mid p\in V_y \text{ and } y\cdot p\in Y\}.
\end{equation*}
As $Y$ is an open subset of $X$ and because the action is continuous, we have that $V'_y$ is an open subset of $V_y$. It is a neighborhood of $e$ satisfying the requirements of Definition \ref{defin:1.6} for the local action on $Y$. We use the argument as in the remark after Proposition \ref{prop:1.8}. 
\ebew

We now discuss two simple examples. More examples that illustrate this procedure are found in \cite{La-VD2, La-VD3}.

\voorb
i) For the group $G$ we take the group with two elements, the identity $e$ and an element $p$ satisfying $p^2=e$.  Let $X$ be the open interval $]-2,2\,[$ in $\mathbb R$. Define a global right action of $G$ by $x\cdot p=-x$ for all $x\in X$. Let $Y$ be the open subset $\,]-1,2\,[$ of $X$. We look at the restriction of the global action to this subset, as defined in \ref{prop:1.7}.
\ssnl
By definition we have $(y,e)\in \Gamma$  and $y\cdot e=y$ for all $y\in Y$. Further $(y,p)\in \Gamma$ if $y$ and $-y$ belong to $Y$. This means that $y\in\,]-1,1[$. For the restricted action we still have $y\cdot p=-y$ for $y\in\,]-1,1[$.
\ssnl
We obtain the following local action of the group $G$ on $Y$ where
\begin{itemize}[noitemsep]
\item[a)] $Y=\,]-1,2\,[$, 
\item[b)]$\Gamma$ is the union of $Y\times\{e\}$ and $Y_0\times \{p\}$ where $Y_0=\,]-1,1[$,
\item[c)] the action is given by $y\cdot p=-y$ for $y\in Y_0$.
\end{itemize}
\ssnl
ii) This action is equivalent with the local right action of $G$ on $\mathbb R$ given by 
\begin{equation*}
x\cdot q=
\begin{cases}
x &\text{ for } q=e \text{ and all } x\in \mathbb R\\
x\inv& \text{ for }q=p \text{ and } x>0.
\end{cases}
\end{equation*}
\evoorb

This example raises the following questions:
\begin{itemize}[noitemsep]
\item[i)] Is every local action of the two element group on $\mathbb R$  the restriction of a global action?
\item[ii] The same question for a every local action on any $X$?
\end{itemize}

The following example shows that the answer to the first question is negative.

\prop\label{prop:1.12}
Consider the local action of the two element group on $\mathbb R$ given by 
\begin{equation*}
x\cdot q=
\begin{cases}
x &\text{ for } q=e \text{ and all } x\in \mathbb R\\
x & \text{ for }q=p \text{ and } x>0.
\end{cases}
\end{equation*}
This local right action can not be realized as the restriction to an open subset of a global action.
\eprop
\bew
Suppose that $X$ is a locally compact space with a global action of $G$. The global action is given by an isomorphism $\alpha$ of $X$.
\ssnl
Assume that $X$ containes an open subset $Y$, homeomorphic to $\mathbb R$. Then $\alpha(Y)$ is another open subset, homeomorphic to $\mathbb R$. Elements that belong to $Y\cap \alpha(Y)$ are elements $y$ of $Y$ that satisfy $y> 0$. This is  because of the assumption that the restriction to $Y$ is the given local action. 
\ssnl
For these elements we have $\alpha(y)=y$. The element $0$ can be approximated by elements $y\in Y$ with $y>0$. Because $\alpha(y)=y$ for all $y>0$ , also $\alpha(0)=0$. So also $0\in Y\cap\alpha(Y)$. As this is an open subset, we must also have negative elements in $Y\cap \alpha(Y)$. This is a contraction. 
\ebew

Clearly this local action is also obtained as the restriction of the trivial action of $G$ on $\mathbb R$, but this is another type of restriction. We restrict the action on $\mathbb R\times G$ to $\Gamma$.
\nl
\bf Other types of restrictions \rm
\nl
Given a local right action on $\Gamma$ in $X\times G$, then we can look for open subsets $\Gamma_0$ of $\Gamma$ and restrict the action to this open subset. It is needed that $X\times \{e\}$ is still contained in $\Gamma_0$. We need that, given a pair $(x,p)\in \Gamma_0$ and any other element $q\in G$, then $(x,pq)\in \Gamma$ if and only if $(x\cdot p, q)\in \Gamma$. We have an (almost trivial) example of this in Proposition \ref{prop:1.12}.
\ssnl
We can also \emph{restrict the action to a closed subgroup}.

\defin\label{defin:1.10}
Assume that we are given a local right action of $G$ on $X$. Let $H$ be a closed subgroup of $G$. Define $\Gamma_H$ as the set of pairs $(x,h)$ in $X\times H$ satisfying $(x,h)\in \Gamma$. Then $(x,h)\mapsto x\cdot h$ is a local right action of $H$ on $X$. We call it the restriction of the action of $G$ to the subgroup $H$.
\edefin

It is straightforward to show that the restriction to a closed subgroup of a proper action is still a proper action. For the notion of locally homeomorphic actions, there is no such result. But if the action of the subgroup is locally homeomorphic, some interesting consequences can be shown for the action of the larger group, see e.g.\ Proposition \ref{prop:6.15a} and some of its consequences. 
\nl
\bf Derived actions \rm
\nl
\oldcomment{\rood Refer to later sections where we treat polynomial functions for these derived actions.
\nl}
Let a local right action of $G$ on $X$ as in Definition \ref{defin:1.1} be given, then we can associate various other actions.
\ssnl
First we have left and right multiplication with $G$ on the second factor in the Cartesian product $X\times G$. They are global actions and can be restricted to the open subset $\Gamma$ as in Proposition \ref{prop:1.7}. We get commuting left and a right local actions of $G$ on $\Gamma$. 
\ssnl
We can also consider the local right action of $G$ on the first factor in $X\times G$ and again we can restrict it to $\Gamma.$ Finally, we can combine two of these actions as follows.

\prop\label{prop:1.14}
There is a  local right action of $G$ on $X\times G$ defined by $(x,p)\tl_1 q=(x\cdot q,q\inv p)$ on the subset of elements $(x,p,q)$ in $(X\times G)\times G$ satisfying $(x,q)\in \Gamma$. The subset $\Gamma$ of $X\times G$ is invariant for this local action.
\eprop

\bew
i) Define a subset $\Gamma_1$ of $(X\times G)\times G$ of elements $(x,p,q)$ satisfying $(x,q)\in \Gamma$. Then $\Gamma_1$ is an open subset of $(X\times G)\times G$ because $\Gamma$ is an open subset of $X\times G$. Further we define  $(x,p)\tl_1 q=(x\cdot q,q\inv p)$  in $X\times G$ for $(x,p,q)\in \Gamma_1$. 
\ssnl
The map $(x,p,q)\mapsto (x\cdot q,q\inv p)$ is a continuous map from $\Gamma_1$ to $X\times G$. We clearly have $(x,p,e)\in \Gamma_1$ and $(x,p)\tl_1 e=(x,p)$. Next suppose that $(x,p,q)\in \Gamma_1$ and take another element $r\in G$. Then $(x,p,qr)\in \Gamma_1$ if and only if $(x,qr)\in\Gamma$. On the other hand, $(x\cdot q,q\inv p,r)\in \Gamma_1$ if and only if $(x\cdot q,r)\in \Gamma$. For the original action we know that $(x,qr)\in\Gamma$ if and only if $(x\cdot q,r)\in \Gamma$. This shows that $(x,p,qr)\in \Gamma_1$ if and only if $((x,p)\tl_1 q,r)\in\Gamma_1$. Also in that case we have
\begin{align*}
(x,p)\tl_1 (qr)
&=(x\cdot (qr),(qr)\inv p)\\
&=((x\cdot q)\cdot r, r\inv q\inv p)\\
&=(x\cdot q,q\inv p)\tl_1 r\\
&=((x,p)\tl_1 q)\tl_1 r.
\end{align*}
Therefore we get a local right action of $G$ on $X\times G$.
\ssnl
ii) Now let $(x,p,q)\in \Gamma_1$. Because $(x,q)\in\Gamma$, by the compatibility condition ii) in Definition \ref{defin:1.1} we have
\begin{equation*}
((x\cdot q),q\inv p)\in \Gamma  \text{ \ if and only if \ } (x, qq\inv p)=(x,p)\in \Gamma.
\end{equation*}
Therefore, if also $(x,p)\in \Gamma$ we have $(x,p)\tl_1 q\in \Gamma$ and the local action of $G$ on $X\times G$ will leave $\Gamma$ invariant. It will yield a local action of $G$ on the subset $\Gamma$ of $X\times G$.
\ebew
We still denote this local action of $G$ on $\Gamma$ with $\tl_1$. As we mentioned already, we can also define another local action $G$ on $\Gamma$. It is defined by the formula $(x,p)\tl_2 q=(x,pq)$ for $x\in X$ and $p,q\in G$. We get the following relation.

\prop\label{prop:1.15}
Consider the global right action of $G$ on $X\times G$ given by $(x,p)\tl_2 q=(x,pq)$. We can restrict it to the open subset $\Gamma$ as in Proposition \ref{prop:1.7} and we get a local right action of $G$ on $\Gamma$. We use again $\tl_2$ for this restriction. The map $\gamma$, defined from $\Gamma$ to $X\times G$ in Proposition \ref{prop:1.4} by $\gamma(x,p)=(x\cdot p,p\inv)$, intertwines  the local action $\tl_1$ from the previous proposition with the local action $\tl_2$ obtained here.
\eprop

\bew
i) It is obvious that $(x,p,q)\mapsto (x,p)\tl_2 q:=(x,pq)$ is a right global action of $G$ on $X\times G$. We can consider its restriction to $\Gamma$ as in Proposition \ref{prop:1.7}.
\ssnl
ii) Consider the map $\gamma$ on $\Gamma$, defined in Proposition \ref{prop:1.4} as $\gamma(x,p)=(x\cdot p,p\inv)$. Recall that it is a homeomorphism from $\Gamma$ to itself.
\ssnl
Given $(x,p)\in \Gamma$ and $(x,q)\in \Gamma$ we find
\begin{align*}
(\gamma(x,p))\tl_2 q
&=(x\cdot p,p\inv)\tl_2 q\\
&=(x\cdot p,p\inv q)\\
&=\gamma(x\cdot q, q\inv p)\\
&=\gamma((x,p)\tl_1 q).
\end{align*}
This proves that $\gamma$ intertwines these two actions.
\ebew

As we have seen in Proposition \ref{prop:1.12} there are partial actions that do not arise as the restriction to an open set of a global action as in Proposition \ref{prop:1.7}. On the other hand,  the two results show that, in a sense, this is \emph{almost} the case.  \oldcomment{\rood Do we have an example?}{}
\ssnl
Indeed, we see that the local action $\tl_1$ of $G$ on $\Gamma$ is equivalent to the restriction to $\Gamma$ of the global action $\tl_2$ of $G$ on $X\times G$. Of course, we cannot conclude from this that the original local right action of $G$ on $X$ is still the restriction of a global action. But nevertheless, this fact shows that, to a certain extend, local actions can be studied by looking at the restrictions of global actions.
\ssnl
It also indicates that local actions obtained from a global action as in Proposition \ref{prop:1.7} can serve as a model for the study of local actions.
\ssnl 
Observe that the global right action of $G$ on $X\times G$ given by $(x,p)\tl_2 q=(x,pq)$ is a proper action. Hence that is also true for the restriction to the open set $\Gamma$ and by the result above, the local action $\tl_1$ of $G$ on $\Gamma$ as defined in the Proposition \ref{prop:1.15}, is again a proper action.
\ssnl
We have a number of generalizations of the previous derived actions.
\nl
\bf Actions associated to admissible pairs \rm
\nl
The purpose of this paper is to provide the basis for the application of local actions and the associated polynomial functions in the theory of bicrossproducts as studied in  \cite{La-VD3}. 

\section{\hspace{-17pt}. Induced actions on functions}\label{s:ind-fun} 

Let $(x,p)\mapsto x\cdot p$ be a local right action of $G$ on $X$, defined on the open subset $\Gamma$ of $X\times G$. In this section we look at the behavior of  functions $f$ in $C_c(X)$ with respect to such a local action. 
\ssnl
In the case of a global action, the situation is obvious. For every element $p\in G$ we have the function $p\cdot f$ defined on $X$ as $x\mapsto f(x\cdot p)$. It is again a function in $C_c(X)$ and the map $G\times C_c(X)\to C_c(X)$ given by $(p,f)\mapsto p\cdot f$ is a global action of $G$ on $C_c(X)$.
\ssnl
In the general case of a local action, the situation is  more complex and there are some important issues to be considered. The following result is crucial for our approach.

\prop\label{prop:2.2}
Let $f\in C_c(X)$ and $p\in G$. Define a function $p\cdot f$ on $X$ by
\begin{equation*}
(p\cdot f)(x)=
\begin{cases}
f(x\cdot p) &\text{ if } (x,p)\in \Gamma\\
0&\text{ otherwise}.
\end{cases}
\end{equation*}
Denote by $C$ the support of $f$. If $(x,p\inv)\in \Gamma$ for all $x\in C$, then $p\cdot f$ is continuous with compact support. In that case, its support is equal to the set of points  $x\in X$ satisfying $(x,p)\in\Gamma$ and $x\cdot p\in C$. 
\eprop
\bew
i) Denote  by $C_p$ the set of points $x\in X$ satisfying $(x,p)\in\Gamma$ and $x\cdot p \in C$. By our assumption, we have $x\in C_p$ if and only if $x=y\cdot p\inv$ for some $y\in C$.  Therefore $C_p$ is a compact subset of $X_p$ as image of the compact set $C$ under the continuous map $y\mapsto y\cdot p\inv$.
\ssnl
ii) If  $(p\cdot f)(x)\neq 0$ we must have $(x,p)\in \Gamma$ and $f(x\cdot p)\neq 0$ so that $x\cdot p\in C$ and $x\in C_p$. We see that $p\cdot f$ is continuous on $X_p$ (as the action is continuous and $f$ is continuous) and on $X\setminus C_p$ (where it is $0$). As $C_p$ is a closed subset of $X_p$ it follows that $p\cdot f$ is continuous on all of $X$. The support of $p\cdot f$ is clearly $C_p$.
\ebew

We will write $C\cdot p\inv$ for the set $C_p$. Recall that, when we write $C\cdot p\inv$, it is understood that this is well-defined, i.e.\ that $(y,p\inv)\in \Gamma$ when $y\in C$.
\ssnl
Although we have  defined $p\cdot f$ for all $p$ in general, we will only use it when all elements $x$ in the support  $C$ of $f$ satisfying $(x,p\inv)\in\Gamma$ so that $p\cdot f\in C_c(X)$.

\opm\label{opm:2.3}
It will sometimes be convenient to use $x\mapsto f(x\cdot p)$ for the function $p\cdot f$ when the support $C$ of $f$ satisfies $(x,p\inv)\in\Gamma$ for all $x\in C$. Remark that anyway, $(p\cdot f)(x)=f(x\cdot p)$ when $x\cdot p$ is defined and that otherwise, $(p\cdot f)(x)=0$. 
\eopm

\defin\label{defin:2.3}
Let $f\in C_c(X)$ with support $C$. Denote by $V_f$ the set of elements $p\in G$ satisfying $(x,p\inv)\in \Gamma$ for all $x\in C$.
\edefin

Recall that for all $p\in V_f$,  now $p\cdot f$ belongs to $C_c(X)$  and that its support is $C\cdot p\inv$.

\prop \label{prop:2.3}
The set $V_f$ is a neighborhood  of the identity in $G$.
\eprop

\bew
i) Because $x\cdot e$ is defined for all $x\in X$, and hence all $x\in C$, we have $e\in V_f$.
\ssnl
ii) To prove that $V_f$ is open, take an element $p_0\in V_f$. For all $x\in C$ we have $(x,p_0\inv)\in \Gamma$. Because $\Gamma$ is an open subset of the Cartesian product, we have a neighborhood  $U_x$ of $x$ in $X$ and a neighborhood $V_x$ of $p_0$ in $G$ so that $(y,p\inv)\in\Gamma$  for all  $y\in U_x$ and all $p\in V_x$. The sets $U_x$ form an open cover of the compact set $C$. Take a finite subcover
\begin{equation*}
C\subseteq U_{x_1}\textstyle\bigcup U_{x_2}\textstyle\bigcup \dots \textstyle\bigcup U_{x_n}
\end{equation*}
and define $V=\bigcap_j V_{x_j}$.
\ssnl
iii) Clearly $V$ is a neighborhood of $p_0$. Now take $p\in V$. For any element $x\in C$ we can choose $j$ so that $x\in U_{x_j}$. Then $p\in V_{x_j}$ as $V\subseteq V_{x_j}$. So $(x,p\inv)\in\Gamma$ and it follows that $p\in V$. This proves that $V$ is open.
\ebew

Keep in mind that, if $(p\cdot f)(x)\neq 0$, then $(x,p)\in \Gamma$ and $(p\cdot f)(x)=f(x\cdot p)$. 

\opm\label{opm:2.5}
Suppose that we start from a local right action of $G$ on $X$ and we restrict it to an open subset $Y$ as in Proposition \ref{prop:1.7}. Let $f$ be a function in $C_c(X)$ and assume that its support $C$ is a subset of $Y$. The set $V_f$  is the set of elements $p\in G$ satsifying $(y,p\inv)\in\Gamma$ for all $y\in C$. For the restricted action, we have a possibly smaller set $V'_f$ of those elements $p\in V_f$ so that also $y\cdot p\inv\in Y$ for all $y$ in $C$. 
\eopm
\oldcomment{\rood Should we not review this part - Think about it!}
\nl
\bf The associated local left action of $G$ on $C_c(X)$ \rm
\nl
We now show that we have a local left action of $G$ on $C_c(X)$ in the following sense.

\prop\label{prop:2.6}
Let $f$ be a function in  $C_c(X)$.  Then $e\in V_f$ and $e\cdot f=f$. Now take any $p\in V_f$ and let $g=p\cdot f$. For $q$ in $G$ we have $q\in V_g$ if and only if $qp\in V_f$. In that case $(qp)\cdot f=q\cdot g$. 
\eprop

\bew
i) The first statement is trivial because $(x,e)\in\Gamma$ and $x\cdot e=x$ for all $x\in X$.
\ssnl
ii) Let $f\in C_c(X)$ with  support denoted as before by $C$. Let $p\in V_f$. Then $p\cdot f$ is well-defined in $C_c(X)$ by definition (cf.\ Proposition \ref{prop:2.2}). The support of $g$ is the set of points $y\cdot p\inv$ where $y\in C$. So $V_g$ consists of elements $q\in G$ with the property that $(y\cdot p\inv,q\inv)\in\Gamma$. This is equal to the set of points $q$ satisfying $(y,(p\inv q\inv))\in\Gamma$. It follows that $q\in V_g$ if and only if $qp\in V_f$ . 
\ssnl
iii) We now prove that $((qp)\cdot f)(x)=(q\cdot g)(x)$ for all $x$ if $p\in V_f$ and $q\in V_g$. We have to consider different possiblities for $x$. First assume that $(x,q)\in\Gamma$ and also $(x\cdot q,p)\in \Gamma$.Then we have
\begin{equation*}
(q\cdot g)(x)=g(x\cdot q)=(p\cdot f)(x\cdot q)=f(x\cdot q\cdot p)=f(x\cdot (qp))=((qp)\cdot f)(x).
\end{equation*}
If $(x,q)\in \Gamma$ and $(x\cdot q,p)\notin\Gamma$ then also $(x,qp)\notin \Gamma$ and  we  have 
\begin{equation*}
(q\cdot g)(x)=g(x\cdot q)=(p\cdot f)(x\cdot q)=0
\end{equation*}
as well as $(qp\cdot f)(x)=0$. 
\ssnl
Finally assume that $(x,q)\notin\Gamma$. Then $(q\cdot g)(x)=0$. Now suppose that $(qp\cdot f)(x)\neq 0$. This means that $x\in C\cdot (qp)\inv$. Now we use that $p\in V_f$ so that $C\cdot p\inv$ is defined. This implies that $x\in (C\cdot p\inv)\cdot q\inv$ and hence $x\cdot q$ is defined. This contradicts the assumption that $(x\cdot q)\notin\Gamma$. So we also have $((qp)\cdot f)(x)=0$. As also $(q\cdot g)(x)=0$, we again get equality.
\ebew

We will need the following weak form of continuity of this action. 

\prop\label{prop:2.7}
For each $x\in X$, the function $v\mapsto (v\cdot f)(x)$ is continuous on $V_f$.
\eprop
\bew
Fix $x\in X$ and $v_0$ in $V_f$.
\ssnl
i) First assume that $(x,v_0)\in \Gamma$. Then we can choose a neighborhood $W$ of $v_0$ so that $(x,w)\in \Gamma$ for all $w\in W$. By taking the intersection with $V_f$, we can assume that $W\subseteq V_f$. For $w\in W$ we have
\begin{equation*}
(w\cdot f)(x)=f(x\cdot w)
\end{equation*}
because $(x,w)\in\Gamma$. Then $w\mapsto f(x\cdot w)$ is continuous on $W$ because the action is continuous and $f$ is continuous. In particular, it is continuous in $v_0$ and this is true for any $v_0$ in $V_f$.
\ssnl
ii) Next suppose that $(x,v_0)\notin \Gamma$. We claim that then $x\notin C\cdot v_0\inv$. Indeed, if $x\in C\cdot v_0\inv$ then it is of the form $y\cdot v_0\inv$ with $y\in C$. This would imply $(x,v_0)\in\Gamma$. Now we use that the set $C\cdot v_0\inv$ is closed so that its complement is open. By the continuity of the action, it follows that there is a neighborhood $U$ of the identity so that $x\cdot U\notin C\cdot v_0\inv$. Then $x\notin C\cdot v_0\inv U\inv$. Let $W=Uv_0$ and let $w\in W$. If $(x,w)\in\Gamma$, then $(w\cdot f)(x)=f(x\cdot w)=0$. If $(x,w)\notin \Gamma$ also $(w\cdot f)(x)=0$. We see that $(w\cdot f)(x)=0$ for all $w$ in a neighborhood $W$ of $v_0$ and $v\mapsto (v\cdot f)(x)$ is continuous. 
\ebew

Later in this paper, we will have to deal with  local right  actions of a group $G$ and the restriction to a closed subgroup $H$. We include the following remark for this case.

\opm\label{opm:2.8}
 Suppose we have a function $f\in C_c(X)$  and consider the neighborhood $V_f$ of the identity. It is clear that the intersection $H\cap V_f$  is
  the neighborhood of the identity associated with $f$ for the restriction of the action to $H$. 
\eopm

We have similar properties for local left actions. In this case we get a local right action on $C_c(X)$. 
\snl
Let us also briefly look at the case of commuting left and right actions as defined in Definition \ref{defin:1.2}. First look at an argument in the proof of  Proposition  \ref{prop:2.3} for this case.

\lem\label{lem:2.9a}
Assume that we have a  local left action of $H$ and a  local right action of $K$ that are commuting in the sense of Definition \ref{defin:1.2}. Given a compact subset $C$ of $X$, we have neighborhoods $U$ of the identity in $H$ and $V$ of the identity in $K$ so that for all $x\in C$, $h\in U$ and $k\in V$ such that both $h\cdot (x\cdot k)$ and $(h\cdot x)\cdot k$ are defined.
\elem

\bew
We can use the argument from the proof of Propostion \ref{prop:2.3} and we find that the set of pairs $(h,k)\in H\times K$ so that $h\cdot (x\cdot k)$ and $(h\cdot x)\cdot k$ are defined for all $x\in C$ is a neighborhood of the identity in the group $H\times K$. Such a neighborhood contains a neighborhood of the form $U\times V$.
\ebew

Using this property, we can associate the commuting local actions of $H$ and $K$ on $C_c(X)$.

\prop\label{prop:2.10}
Let $f$ be a function in $C_c(X)$.  Denote its support by $C$. Let $U$ and $V$ be as in the lemma above. Let $h\in U\inv$ and $k\in V\inv$. 
Define a function  $g$ on $X$ by
\begin{equation*}
g(x)=
\begin{cases} f(h\cdot x\cdot k) & \text{ if } h\cdot x\cdot k \text{ is defined } \\
0 & \text{ otherwise.} \end{cases}
\end{equation*}
Then $g\in C_c(X)$.
\eprop
\bew
The argument is the same as  in the proof of Proposition \ref{prop:2.2}. On the set $X_{h,k}$ of points $x$ where $h\cdot x\cdot k$ is defined, we have continuity of $g$ because $f$ is continuous and because the actions are continuous. On the other, hand when $x$ is not an element in $h\inv\cdot C\cdot k\inv$ we have $g(x)=0$.  We see that the function is continuous and has compact support in the open set $X_{h,k}$. Therefore it is continuous on all of $X$.
\ebew
We denote this function as $k\cdot f\cdot h$. 
\ssnl
As we have seen in Proposition \ref{prop:1.3}, a pair of commuting left and right local actions gives rise to a single right local action. Then the function defined in Proposition \ref{prop:2.10} is just as the one defined in Proposition \ref{prop:2.2}.

\section{\hspace{-17pt}. Polynomial functions}\label{s:pol} 

In this section, we initiate the study of polynomial functions. We consider a  local right action of the group $G$ on the space $X$ as in definition \ref{defin:1.1}. We define the notion for functions in $C_c(X)$. 
We use the notations of the previous sections.

\defin\label{defin:3.1}
Let $f$ be a function in $C_c(X)$. If there is a neighborhood  $V$ of the identity in $G$, contained in $V_f$, and a \emph{finite-dimensional} subspace $F$ of $C_c(X)$ so that  the functions $v\cdot f$ with $v\in V$ all belong to $F$, then we call $f$ \emph{polynomial}. 
\edefin

It is clear that a scalar multiple of a polynomial function, as well as the sum and the product of two polynomial functions is again polynomial. 

\notat\label{notat:3.2} We use $\mathcal P(X)$ for the space of polynomial functions in $C_c(X)$. When we take for $X$ the group $G$ and  right multiplication for the (global) action, we use $\mathcal P(G)$ for the set of polynomial functions in $C_c(G)$.
\enotat

Strictly speaking we should include the symbol $G$ and the action in the notation $\mathcal P(X)$. We omit that for simplicity. We just keep in mind that the notion depends on the local action of the locally compact group $G$. 
\snl
The following is quite obvious, but important for results to come.

\prop\label{prop:3.3}
Suppose that $H$ is a closed subgroup of $G$. If $f$ is a polynomial function in $C_c(X)$ for the action of $G$, it is also polynomial for the restriction of the action to $H$.
\eprop

\bew
Assume $f$ polynomial and take a neighborhood of the identity $V$, contained in $V_f$ and a finite-dimensional subspace $F$ of $C_c(X)$ containing the functions $v\cdot f$ for all $v\in V$. Then $V\cap H$ will be a neighborhood of the identity in $H$ and it is contained in $V_f\cap H$. We have seen in Remark \ref{opm:2.8} that $V_f\cap H$ is the set $V_f$ for the restricted action.  We can further take the same finite-dimensional subspace $F$. Then $w\cdot f\in F$ for all $w\in V\cap H$ as this is a subset of $V$.
\ebew

Later we will prove that, under certain special conditions, it is true that $f$ is polynomial for the action of $G$ when it is polynomial for the action of $H$, see Proposition \ref{prop:6.15a}. 
\ssnl
In Section \ref{s:misc} we will also consider polynomial functions for the derived actions we have in Section \ref{s:loc-act}, see e.g.\ Proposition \ref{prop:6.6}.

\opm\label{opm:3.4}
i) For a local action of a discrete group, all elements in $C_c(X)$ are polynomial functions. This is trivial as we can take for $V$ the singleton $\{e\}$ and for $F$ the space spanned by $f$.  
\ssnl
ii) Also if the group action is trivial, i.e.\ if $\Gamma=X\times G$ and $x\cdot p=x$ for all $x\in X$ and $p\in G$, all elements in $C_c(X)$ are polynomial.  Now we can take $G$ for $V$ and for $F$ again the space spanned by $f$.
\ssnl
iii) When the discrete group $G$ acts on itself by multiplication, we find that $\mathcal P(G)$ is the space of all complex functions with finite support. This follows from Remark i) above. 
\eopm
These are trivial statements but they show that \emph{in general}, the existence of polynomial functions puts no restrictions on neither $X$ or $G$, only on the action of $G$ on $X$.  
\nl
\bf First properties of polynomial functions \rm
\nl
Some of the properties of polynomial functions will be proven by using the following general result.

\lem\label{lem:3.5}
Assume that we have a finite-dimensional subspace $F$ of $C_c(X)$ and that $(f_i)$ is a basis for $F$. Then there exist finitely many complex numbers $c_{jk}$ and elements $x_{jk}$ in $X$ so that 
\begin{equation}
\sum_j c_{jk}f_i(x_{jk})=\delta_{ik}\label{eqn:3.1}
\end{equation}
for all $i,k$. Here $\delta_{ik}$ is the Kronecker $\delta$.
\elem

\bew
Let $n$ be the dimension of $F$. 
For each $x\in X$ we consider the element $(f_i(x))_i$ in $\mathbb C^n$. The span  of these elements must be all of $\mathbb C^n$. Indeed, if  we have numbers $c_i$ so that $\sum_i c_i f_i(x)=0$ for all $X$, these numbers must all be $0$ because the functions $f_i$ are linearly independent. Then for any element $\xi\in \mathbb C^n$ we find numbers $c_j$ and elements $x_j$ in $X$ so that 
\begin{equation*}
\xi_i=\sum_j c_j f_i(x_j)
\end{equation*}
for all $i$. Now take an index $k$ and  let $\xi_i=\delta_{ik}$ for all $i$. Then we get  complex numbers $c_{jk}$ and elements $x_{jk}$ in $X$ so that Equation (\ref{eqn:3.1}) is satisfied.
\ebew

We will  apply this result not only for a finite number of functions in $C_c(X)$ but also for functions on a neighborhood $V$ of the identity in $G$. We formulate the result in the following proposition. 

\prop\label{prop:3.6}
Let $\{\varphi_i\}$ be finitely many functions in $C_c(G)$. Assume that $V$ is  a neighborhood of the identity in $G$ and that the restrictions ${\varphi_i}|_V$ are linearly independent. Then  there exist complex numbers $c_{jk}$ and elements $v_{jk}$ in $V$ so that 
\begin{equation}
\sum_j c_{jk}\varphi_i(v_{jk})=\delta_{ik}\label{eqn:3.2}
\end{equation}
for all $i,k$. 
\eprop

This follows from the lemma or can be proven directly as in the proof of the lemma.
\ssnl

We will now obtain some  properties that characterize polynomial functions in a more convenient way. To do this, we need the following basic property of polynomial functions.

\prop\label{prop:3.7} 
Let $f$ be a polynomial function in $C_c(X)$. If  $p\in V_f$, then $p\cdot  f$ will be again a polynomial function in $C_c(X)$.  
\eprop

\bew  
i) Let $f$ be  a polynomial function in $C_c(X)$, fix $p\in V_f$ and take  $g=p\cdot f$. As in Definition \ref{defin:3.1},  let $V$ be a neighborhood of the identity, contained in $V_f$, and $F$ a finite-dimensional subspace of $C_c(X)$ containing all the functions $v\cdot f$ with $v\in V$. 
\ssnl
ii) Because $p\in V_f$ and $V_f$ is open, we have a neighborhood $V'$ of the identity so that $V'p\subseteq V_f$. By taking it small enough, we can assume that also $V'\subseteq pVp\inv$.
\ssnl
iii) Take $v\in V'$. Because $p$ and $vp$ are both in $V_f$ we know by Proposition \ref{prop:2.6} that $v\in V_g$ and $v\cdot g=(vp)\cdot f$. Further $vp=p(p\inv vp)$ and both $vp$ and $p\inv vp$ are elements in $V_f$. If we take $h=(p\inv vp)\cdot f$, again by Proposition \ref{prop:2.6},  we also have that $p\in V_h$ and 
\begin{equation*}
v\cdot g=(vp)\cdot f=p\cdot h.
\end{equation*}
All the functions $(p\inv vp)\cdot f$, with $v\in V'$ belong to a finite-dimensional subspace $F_1$ of $C_c(X)$, contained in $F$. Moreover, $p\in V_h$ for all $h$ in this subspace. Denote by $F'$ the space of functions $p\cdot h$ where $h\in F_1$. Then $F'$ is a finite-dimensional subspace of $C_c(X)$ containing all functions $v\cdot g$ with $v\in V'$. This proves that $g$ is polynomial.

\ebew

The proof in the global case is much easier. We simply replace $V$ by $pVp\inv$ and $F$ by the space of functions $p\cdot h$ where $h\in F$. 
\ssnl
The following property is essentially a reformulation of the definition. 

\prop\label{prop:3.8}
Let $f$ be a polynomial function in $C_c(X)$. Then there is a neighborhood $V$ of the identity in  $G$ contained in $V_f$ and finitely many  functions $\varphi_i\in C_c(G)$  and $f_i\in C_c(X)$ so that 
\begin{equation}
v\cdot f=\sum_i\varphi_i(v) f_i\label{eqn:3.3}
\end{equation}
for all 
$v\in V$. The functions $f_i$ can be chosen to be polynomial.
\eprop

\bew
i) Suppose that $f$ is a non-zero polynomial function in $C_c(X)$. Choose a neighborhood $W$ of the identity in $G$ contained in $V_f$ and a finite-dimensional subspace $F$ of $C_c(X)$ containing the functions $w\cdot f$ for all $w\in W$. Take a basis $(f_i)$ for $F$ and define complex valued functions $\psi_i$ on $W$ by
\begin{equation*}
w\cdot f=\sum_i \psi_i(w) f_i
\end{equation*}
for all  $w\in W$.
Because the functions $f_i$ are linearly independent, by Lemma \ref{lem:3.5} we find numbers $c_{jk}$ and elements $x_{jk}\in X$ so that $\sum_j c_{jk}f_i(x_{jk})=\delta_{ik}$. 
Then
\begin{equation*}
\psi_k(w)=\sum_{i,j} c_{jk}\psi_i(w)f_i(x_{jk})=\sum_j c_{jk}(w\cdot f)(x_{jk}).
\end{equation*}
We have shown that the functions $w\mapsto (w\cdot f)(x)$ are continuous on $V_f$ for all $x$, see Proposition \ref{prop:2.7}. It follows that $\psi_k$ is continuous on $W$ for all $k$.
 \ssnl
ii) Now choose a neighborhood $V$ of the identity with compact closure contained in $W$ and a function $\varphi$ in $C_c(G)$ with support in $W$ and so that $\varphi(v)=1$ when $v\in V$. Define functions $\varphi_i$ on all of $G$ by $\varphi_i(w)=\varphi(w)\psi_i(w)$ for $w\in W$ and $\varphi_i(w)=0$ for $w\notin W$. These functions are in $C_c(G)$  because $\psi_i$ is defined and continuous on $W$ and $\varphi$ has compact support in $W$. Furthermore 
\begin{equation}
v\cdot f=\sum_i \varphi_i(v)  f_i \label{eqn:3.4}
\end{equation}
for all $v\in V$ because $\varphi(v)=1$ for $v\in V$. 
\ssnl
iii) Next assume that the functions $\varphi_i$ are linearly independent on $V$. By Proposition \ref{prop:3.6} we have  numbers $c_{jk}$ and elements $v_{jk}\in V$ so that 
$ \sum_j c_{jk} \varphi_i(v_{jk})=\delta_{ik}$
for all $i,k$. Then
\begin{equation*}
\sum_j c_{jk}(v_{jk}\cdot f)=\sum_{i,j} c_{jk}\varphi_i(v_{jk})f_i=f_k
\end{equation*}
for all $k$. 
By Proposition \ref{prop:3.7} we know that $v_{jk}\cdot f$ is again polynomial and therefore also $f_k$ is polynomial for all $k$. 
\ebew

We will need the following relation between $V_f$ and the sets $V_{f_i}$.

\prop\label{prop:3.9}
Consider the result of the previous proposition. Assume that the functions $(f_i)$ are linearly independent. If  $p\in G$ and  $pV\subseteq V_f$, then $p\in V_{f_i}$ for all $i$. 
\eprop

\bew
Let $p\in G$. Suppose that  $pv\in V_f$ for all $v\in V$.
\ssnl
Let $v\in V$. Because $pv\in V_f$ we have that $C_f\cdot (pv)\inv$ is defined where $C_f$ is the support of $f$. As $v\in V_f$ we know that also $C_f\cdot v\inv$ is defined. Then we get that $C_f\cdot v\inv\cdot p\inv$ is defined and because $C_f\cdot v\inv=C_{v\cdot f}$ we get that $C_{v\cdot f}\cdot p\inv$ is defined. This holds for all $v\in V$.
\ssnl
Now let $x\in C_{f_k}$. We claim that there is a $v\in V$ so that $x\in C_{v\cdot f}$. This will imply that $C_{f_k}\cdot p\inv$ is defined and hence $p\in V_{f_k}$. This will yield the result.
\ssnl
To prove the claim, we use that there are elements $v_{jk}$ in $V$ and numbers $c_{jk}$ such that 
\begin{equation*}
f_k=\sum_j c_{jk}(v_{jk}\cdot f). 
\end{equation*}
This was shown in the proof of the previous proposition. Therefore, if $f_k(x)\neq 0$ there must be some $v$ so that $(v\cdot f)(x)\neq 0$. In other words, if $x\in C_{f_k}$ there is a $v\in V$ such that $x\in C_{v\cdot f}$. This proves the claim. 
\ebew

Observe that the function $f$ belongs to the space spanned by the functions $f_i$. Indeed, just take $v=e$ in Equation (\ref{eqn:3.4}). 
\ssnl
On the other hand, if $f$ is a function in $C_c(X)$ so that  there is a neighborhood $V$ of the identity and functions $\varphi_i$ and $f_i$ as in the formulation of the proposition, then $f$ is clearly polynomial. So the property is indeed an equivalent formulation of a polynomial function.
\nl
We will often encounter situations of the following kind and the following problem.

\opm\label{opm:3.10a}
i) Assume that we have a local right action of a group $G$ on a space $X$. Let $Y$ be another locally compact space and consider the Cartesian product $X\times Y$. There is the obvious local right action of $G$ on $X\times Y$ given by $(x,y)\cdot p=(x\cdot p, y)$ when $x\cdot p$ is defined. Let $f$ be a polynomial function in $C_c(X\times Y)$ for this action. For each $y\in Y$, the function $x\mapsto f(x,y)$ will be a polynomial function on $X$ for the original action. Indeed, we can write 
\begin{equation*}
f(x\cdot v,y)=\sum_i\varphi_i(v)f_i(x,y)
\end{equation*}
for $v$ in a neighborhood of the identity in $G$, with functions $\varphi_i\in C_c(G)$ and $f_i\in C_c(X\times Y)$. This holds for all $x$ and all $y$. So obviously $x\mapsto f(x,y)$ is polynomial on $X$ for each $y$.
\ssnl
ii) On the other hand, in general, one can not conclude that $f$ is polynomial if that is the case for each of the functions $x\mapsto f(x,y)$. If this is assumed, given $y\in Y$, we can write
\begin{equation*}
f(x\cdot p,y)=\sum_i \varphi_i(v)f_i(x,y)
\end{equation*}
for a neighborhood $V$ of the identity and functions $\varphi_i$ and $f_i$. But the neighborhood depends on $y$, as well as the functions $\varphi_i$. One can not expect in general that there is a single neighborhood contained in all these neighborhoods, nor that one can find functions $\varphi_i$ independent of the choice of $y$.
\ssnl
iii) Only in very special cases, one is able to get the converse result. If we have a function on $X\times Y$ of the form $(x,y)\mapsto f(x)g(y)$ it will be possible. Indeed, if in this case the function is polynomial for each $y$, then $f$ will be polynomial on $X$ and it follows that the function $(x,y)\mapsto f(x)g(y)$ is polynomial for the action on $X\times Y$.
\eopm

This is an important remark and we will have to consider this problem on various occasions. 

\nl
\bf Existence of an open subgroup \rm
\nl

Consider again the result of  Proposition \ref{prop:3.8}. We have written $v\cdot f$ as a sum $\sum_i \varphi_i(v)f_i$ with all the $f_i$ still polynomial. Next we obtain a stronger version of this property.

\prop\label{prop:3.11}
Let $f$ be a non-zero polynomial function in $C_c(X)$. There exists  finitely many polynomial functions $f_i\in C_c(X)$ so that $f$ belongs to the space spanned by these functions, a neighborhood $V$ of $e$ contained in $V_{f_j}$ for all $j$ and functions $\psi_{ij}$ in $C_c(G)$ so that 
\begin{equation*}
v\cdot f_j=\sum_i  \psi_{ij}(v)f_i
\end{equation*}
for all  $v\in V$ and for all $j$.
\eprop

\bew
i) By the Proposition \ref{prop:3.9} we have a neighborhood $V$ of the identity contained in $V_f$, and finitely many functions $\varphi_i\in C_c(G)$ and $f_i\in C_c(X)$ so that $v\cdot f=\sum_i \varphi_i(v)f_i$ for all $v\in V$. We can assume that the functions $\varphi_i$ are linearly independent on $V$. 
As before, choose numbers $c_{jk}$ and elements $v_{jk}\in V$ so that $ \sum_j c_{jk} \varphi_i(v_{jk})=\delta_{ik}$. We have seen in the proof op Proposition \ref{prop:3.8} that $f_k=\sum_j c_{jk}(v_{jk}\cdot f)$ for all $k$.
\ssnl
ii) Given the elements $v_{jk}$ in $V$, we can find a neighborhood $W$ of the identity so that still $Wv_{jk}\in V$ for all these elements $v_{jk}$. By taking it small enough we can assume that $W\subseteq V_{f_k}$ for all $k$. Then for all $k,w\in W$ we have
\begin{equation*}
w\cdot f_k = \sum_j c_{jk}(w\cdot (v_{jk}\cdot f))
= \sum_j c_{jk} ((wv_{jk})\cdot f)
=\sum_{i,j} c_{jk}\varphi_i(wv_{jk})f_i.
\end{equation*}
Now we define the functions $\psi_{ik}$ on $G$ by
$\psi_{ik}(p)=\sum_j c_{jk}\varphi_i(pv_{jk})$
and we get the desired result.
\ssnl
iii) Remember that $f$ belongs to the space spanned by the functions $f_i$ because $f=e\cdot f=\sum_i \varphi_i(e)f_i$.
\ebew

\keepcomment{\rood Should we reconsider this proof in view of the new Proposition \ref{prop:3.9}? At least refer to it? I think this is now ok.
\ssnl}

The main data are in fact the functions $f_i$ in $C_c(X)$, the neighborhood $V$ and the functions $\psi_{ij}$ in $C_c(G)$. We explain this in the following remark. 

\opm
Suppose that we have finitely many functions $f_i$ in $C_c(X)$ and a neighborhood $V$ contained in all $V_{f_i}$ and that there exists functions $\psi_{ij}$ in $C_c(G)$ with the property that
\begin{equation*}
w\cdot f_j=\sum_j \psi_{ij}(w)f_i
\end{equation*}
for all $w\in V$. Let $f$ be of the form $\sum_i \lambda_i f_i$. We know that $f$ is again  polynomial. More precisely, we have $V\subseteq V_f$ and for all $v\in V$
\begin{equation*}
v\cdot f=\sum_j \lambda_j v\cdot f_j=\sum_j \varphi_i(v)f_i
\end{equation*}
where $\varphi_i(v)=\sum_j \lambda_j \psi_{ij}(v)$. 
\ssnl
In Proposition \ref{prop:3.8} we have shown that all polynomial functions are obtained in that way.
\eopm

As we have mentioned already in Remark \ref{opm:3.4}, in general we can not say anything about the space $X$ and the group $G$ if polynomial functions exists. On the other hand, when polynomial functions exist, we can say something about the action.
\ssnl
To begin with, we have the following general result.

\prop\label{prop:3.13}
Assume that $f$ is a non-zero polynomial function in $C_c(X)$. Then there is an open subgroup $G_0$ of $G$ and a finite-dimensional subspace $F$ of $C_c(X)$ so that $q\in V_f$  and $q\cdot f\in F$   for all $q\in G_0$.
\eprop

\bew
i) As in Proposition \ref{prop:3.8}, take a neighborhood $V$ of the identity in $G$,   functions $f_i\in C_c(X)$ and $\psi_{ij}\in C_c(G)$ satisfying
\begin{equation*}
v\cdot f_j=\sum_i \psi_{ij}(v)f_i
\end{equation*}
for  all $v\in V$. Recall also that $f$ belongs to the subspace $F=\text{span}\{f_i\}$. 
\ssnl
ii) The space $F$ is invariant under the action of elements from $V$.  Then it will be invariant under the action of any  $q\in G_0:=\cup_1^\infty\, V^n$. This means that $q\in V_{f_i}$ and $q\cdot f_i\in F$ for all such elements $q$. 
By assuming that $V$ is symmetric, we see that  $G_0$ is an open subgroup. It has the property  that $q\in V_{f_i}$ and the functions $q\cdot f_i$ are  in $F$ for all $q\in G_0$ and all $i$. In particular, also $q\in V_f$ and $q\cdot f\in F$ for all $q\in G_0$.  

\ebew 

\opm\label{opm:3.14}
If we choose a basis $\{f_i\}$ for $F$, we can write $q\cdot f_j=\sum_i \pi_{ij}(q)f_i$ for all $q\in G_0$ and all $j$ where all the functions $\pi_{ij}\in C_c(G_0)$. 
\keepcomment{\rood Explain! Refer to Remark \ref{opm:2.3} in the previous section (not relevant!) and Notation \ref{notat:4.3} in the next section.}{}
\eopm

Remember that an open subgroup is automatically closed and contains the connected component of the identity. So if $G$ happens to be connected, $G_0=G$. 
\ssnl
 
 We can summarize and obtain the main result of this item.
 
\stel\label{stel:3.15}
Consider a right local action of $G$ on $X$. Let $f$ be a non-zero polynomial function in $C_c(X)$. Then there is a compact subset $C$ of $X$ and an open subgroup $G_0$ of $G$ so that $C\cdot G_0=C$ and so that $f$ belongs to a finite-dimensional subspace of $\mathcal F_C$ of functions with compact support in $C$, invariant under the action of $G_0$.
\estel

\bew
From Proposition \ref{prop:3.8} we find a finite-dimensional subspace $F$ of $C_c(X)$ and a neighborhood $V$ of the identity in $G$ with the following properties.
\begin{itemize}[noitemsep]
\item The function $f$ belongs to $F$.
\item The neighborhood $V$ is contained in $V_g$ for all $g\in F$.
\item For all $v\in V$ and $g\in F$ we have that $v\cdot g$ again belongs to $F$.
\end{itemize}
We can assume that $V$ is symmetric and take $G_0=\cup V^n$.
This is an open subgroup $G_0$  and we still have that $G_0\subseteq V_g$ and $q\cdot g\in F$ for all $g\in F$. 
\ssnl
Now let $C$ be the union of the supports of the functions in $F$. For all $x\in C$ and $q\in G_0$ we will have $(x,q)\in \Gamma$ and $x\cdot q\in C$. In other words $C\cdot G_0\subseteq C$. Because $G_0$ is a group we have in fact $C\cdot G_0=C$. All the functions $g\in F$ have support in $C$ so that $F$ is a finite-dimensional subspace of $\mathcal F$.
\ebew

There is an obvious converse result.

\prop\label{prop:3.16}
Given a compact subset $C$ and an open subgroup $G_0$ satisfying \hfill

$C\cdot G_0=C$. Then every function $f\in \mathcal F_C$  in a finite-dimensional $G_0$-invariant subspace will be polynomial.
\eprop

The proof is obvious.
\ssnl

The following remark has a followup in the next section.

\opm\label{opm:3.17}
i) If the group $G$ is connected, the subgroup $G_0$ will be all of $G$. So it will be the same for all polynomial functions. If there are enough polynomial functions, this will probably imply that the action has to be global. However, this is not the case we are interested in. See a remark further in Section \ref{s:concl}. \keepcomment{\rood Don't forget!}
 \ssnl
ii) Suppose now that we have a finite number $(f_i)$ of polynomial functions. Denote by $C_i$ and $G_i$ the related data. If we let $C=\cup C_i$
 and $G_0=\cap G_i$
 we get a compact set $C$ and an open subgroup $G_0$ with $C\cdot G_0=C$.
  Now every $f_i$ will belong to a finite-dimensional $G_0$-invariant subspace of functions in $C_c(X)$ with compact support in $C$.
\eopm

\nl
{\bf Proper and locally homeomorphic actions}
\nl
We now obtain some consequences in the case of these extra proporties of the action.

\prop\label{prop:3.18}
Assume that the action is proper and that $f$ is a non-zero polynomial function in $C_c(X)$. Then the group $G_0$  in Theorem \ref{stel:3.15} is a compact open subgroup of $G$.
\eprop

\bew
We use the notations of that theorem.  From  $C\cdot G_0=C$ we get $(x\cdot q,x)\in C\times C$ for all $x\in C$ and $q\in G_0$. As the map $(x,p)\mapsto (x\cdot p,x)$ is assumed to be proper, we must have a compact subset of $\Gamma$ containing all elements $(x,q)$ for $x\in C$ and $q\in G_0$. This implies that $G_0$ is compact.
\ebew

This result has the following important consequence.

\gev \label{gev:3.19a} 
Non-trivial polynomial functions in $C_c(G)$ for right multiplication as the action can only exist when $G$ has a compact open subgroup.
\egev

Indeed, the right multiplication of $G$ on itself is a proper action as we mentioned after Definition \ref{defin:1.5}. Compare this argument with the one given in Corollary 1.2 in the original paper \cite{La-VD1}.
\oldcomment{This is a different argument as the one given in our original paper.}
\ssnl
One may wonder if the existence of a compact open subgroup is sufficient to guarantee the existence of polynomial functions. This is the case for the action of the group on itself. But it is not true in general. See  Remark \ref{opm:4.13} and Example \ref{voorb:4.15} in the next section.
\nl
Finally, for a locally homeomorphic action, we get the following.

\prop
Assume that we have a locally homeomorphic action. Let $C$ and $G_0$ be as in Theorem \ref{stel:3.15}. Then $C$ is a compact open subset of $X$.
\eprop
\bew
For all $x\in C$ we have a neighborhood of the identity $V_x$ so that $x\cdot V_x$ is defined and open in $X$. Let $W_x$ be the intersection $G_0\cap V_x$. This is an open subset of $X$, containing $x$. Moreover, when $x\in C$ it will be contained in $C$. It follows that $C$ is open.
\ebew

As a consequence we get the following result for the action of $G$ on itself.

\gev
Assume that $f$ is a non-zero polynomial function on $G$ for right multiplication. Then there is a compact open subset $C$, containing the support of $f$, and a compact open subgroup $G_0$ so that $C\cdot G_0=C$.
\egev

It makes sense to assume the existence of a compact open subgroup, without requiring that the action is proper. This would be a weaker condition and is in fact what we do in the next section.
\nl

\keepcomment{\rood Check the iffalse in the tex file!
\nl}

\section{\hspace{-17pt}. Actions of groups with a compact open subgroup}\label{s:cpt-open} 

In the previous section we have seen that for any polynomial function $f$ in $C_c(X)$ there exist a compact subset $C$ of $X$, containing the support of $f$, an open
subgroup $G_0$ so that $C\cdot G_0$ is defined and equal to $C$, and a finite-dimensional subspace of functions with compact support in $C$, invariant under the action of $G_0$ and containing $f$. See Theorem \ref{stel:3.15}. 
\ssnl
We also have proven  that for a proper action, we can further assume that $G_0$ is  a compact open subgroup, see Proposition \ref{prop:3.18}. 

\voorw\label{voorw:4.1}
In the rest of the section, we assume that $G$ has a compact open subgroup.
\evoorw

It is fulfilled when there are non-zero polynomial functions and when the action is proper, but we will not assume that.
\ssnl
Recall also the following convention.

\notat
Let $G_0$ be a compact open subgroup of $G$ . Suppose further that $C$ is a compact subset of $X$ such that $(x,u)\in\Gamma$ and $x\cdot u\in C$ for all $x\in C$ and $u\in G_0$. As already mentioned, we simply write this condition as $C\cdot G_0\subseteq C$.  In fact, because $G_0$ is a group, we will have $C\cdot G_0=C$. 
\enotat

Given a compact open subgroup $G_0$ of $G$ and a compact subset $C$ of $X$ satisfying $C\cdot G_0=C$, we introduce the following objects.

\notat\label{notat:4.3}
We consider the subspace $\mathcal F_C$ of functions in $C_c(X)$ with support contained in $C$. Further we let $\pi_C$ be the representation of $G_0$ on the space $\mathcal F_C$ given by $(\pi_C(u)f)(x)=f(x\cdot u)$ for $x\in C$ and $u\in G_0$. In fact, as we will mostly work with one such a compact subset $C$, we will drop the indices and use $\mathcal F$ and $\pi$.
\enotat

By the representation theory of compact groups, we know that $\pi$ splits into finite-dimen\-sional representations. Every function in such a finite-dimensional component is a polynomial function. And because we assume that $G$ has a compact open subgroup, all polynomial functions arise in that way. Again see Proposition \ref{prop:3.18} and the preceding results in the previous section.
\ssnl
We now consider this property again and discuss some more details under the Condition \ref{voorw:4.1}. We begin with introducing the following notation.

\notat\label{notat:4.4}
Assume that $\gamma$ is an irreducible representation of $G_0$ on a finite-dimensional space with a basis $(\xi_i)_{i=1}^n$. The matrix elements are continuous functions $\gamma_{ij}$ on $G_0$ defined by $\gamma(u)\xi_j=\sum_i \gamma_{ij}(u)\xi_i$ for all $u$ and all $i$. Given a function $f\in\mathcal F$  we can now define functions $f_{ij}$ in $\mathcal F$ by
\begin{equation*}
f_{ij}(x)=\int_{G_0} f(x\cdot u)\gamma_{ij}(u\inv)\,du
\end{equation*}
where we integrate over the normalized Haar measure on $G_0$.
\enotat 

Later it will become clear why we use $u\inv$ in this formula.
\opm
We have to be a bit more careful with the notations. Because $f$ has support in $C$ and $(x,u)\in \Gamma$ for all $x\in C$ and $u\in G_0$, we have $G_0\subseteq V_f$. The function $u\cdot f$ is defined in $C_c(X)$ and will again have its support in $C$. For $x\in C$ we have $(u\cdot f)(x)=f(x\cdot u)$. And for $x\notin C$ we have $(u\cdot f)(x)=0$. Therefore, we can use the formula as above for $x\in C$, while for $x\notin C$ we get $0$ anyway. 
\eopm
We keep this remark in mind  when  we use $x\mapsto f(x\cdot p)$ for the function $p\cdot f$. See also some of the remarks  (e.g.\ Remark \ref{opm:2.3}) we have made in Section \ref{s:ind-fun} about this notation. 

\prop
The functions $f_{ij}$ 
are polynomial and
\begin{equation*}
(u\cdot f_{ij})(x)=\sum_k f_{ik}(x)\gamma_{kj}(u)
\end{equation*}
for all $x\in X$ and all $u\in G_0$. 
\eprop
\bew
We can assume $x\in C$ because when $x\notin C$, we have $0$ left and right.
\ssnl
Then, for all $i,j$ and all $v\in G_0$ we have
\begin{align*}
f_{ij}(x\cdot v)
&=\int_{G_0} f(x\cdot vu)\gamma_{ij}(u\inv)\,du\\
&=\int_{G_0} f(x\cdot u)\gamma_{ij}(u\inv v)\,du\\
&=\int_{G_0} f(x\cdot u)\sum_k\gamma_{ik}(u\inv)\gamma_{kj}(v)\,du\\
&=\sum_k f_{ik}(x)\gamma_{kj}(v).
\end{align*}
\vspace*{-20pt}

\ebew
 
 We see that for each $i$, if they are non-zero, the functions $f_{ij}$ with $j=1,\dots, n$ span a $G_0$-invariant subspace.  The restriction of the representation $\pi$ of $G_0$ to this subspace is a direct sum of components all equivalent to $\gamma$. We can have at most $n^2$  components of this type. This follows from the representation theory of compact groups.
\ssnl
In fact, instead of the matrix elements of $\gamma$, we can use the character $\chi_\gamma$ of $\gamma$. Recall that the character $\chi_\gamma$ is defined as the function $u\mapsto \sum_i\gamma_{ii}(u)$.

\notat
Assume that $\gamma$ is an irreducible representation of $G_0$. Let $\chi_\gamma$ denote its character. Then we define $P_\gamma f$ for $f\in\mathcal F$ by
\begin{equation*}
(P_\gamma f)(x)=\frac1{n_\gamma}\int_{G_0} f(x\cdot u)\chi_\gamma(u\inv)\,du
\end{equation*}
where $n_\gamma$ is the dimension of the representation $\gamma$. The function $P_\gamma f$ again belongs to $\mathcal F$. It is a polynomial function on $X$.
\enotat

If $\gamma$ acts on a vector space with a basis $(\xi)_i$ and if $\gamma_{ij}$ are the corresponding matrix elements, then $P_\gamma f=\frac1{n_\gamma}\sum_i f_{ii}$ where the functions $f_{ii}$ are defined as in 
Notation \ref{notat:4.4}. We clearly have 
\begin{equation*}
(P_\gamma f)(x\cdot u)=\frac1{n_\gamma}\sum_{ij} f_{ij}(x)\gamma_{ji}(u)
\end{equation*}
for all $x$ and all $u\in G_0$. 

\prop
For $f$ in $\mathcal F$ and $\gamma$ we have $P_\gamma (P_\gamma f)=P_\gamma f$. On the other hand, if $\gamma'$ is another irreducible representation of $G_0$, not equivalent to $\gamma$, then $P_{\gamma'}(P_{\gamma} f)=0$. 
\eprop
\bew
i) Given $\gamma$ and $f$ we have
\begin{align*}
(P_\gamma P_\gamma f)(x)
&=\frac1{n_\gamma}\int_{G_0} (P_\gamma f)(x\cdot u)\chi_\gamma(u\inv)\,du\\
&=\frac1{n_\gamma^2}\sum_{ij}\int_{G_0}  f_{ij}(x)\gamma_{ji}(u)\chi_\gamma(u\inv)\,du.
\end{align*}
Now
\begin{equation*}
\int_{G_0} \gamma_{ji}(u)\chi_\gamma(u\inv)\,du
=\sum_k \int_{G_0} \gamma_{ji}(u)\gamma_{kk}(u\inv)\,du.
\end{equation*}
By the orthogonality relations of matrix elements, this will be equal to $0$ except if $i=j$. Then it will be equal to $n_\gamma$. Therefore
\begin{equation*}
(P_\gamma P_\gamma f)(x)=\frac1{n_\gamma}\sum_i f_{ii}(x)=P_\gamma(x).
\end{equation*}
\ssnl
ii) Now assume that we have $\gamma'$, not equivalent with $\gamma$. Then $P_{\gamma'}P_\gamma(f)=0$ because
\begin{equation*}
\int_{G_0} \gamma_{ji}(u)\chi_{\gamma'}(u\inv)\,du=0,
\end{equation*}
again by the orthogonality relations of the matrix elements.
\ebew

If on the other hand, $f$ belongs to a finite-dimensional invariant subspace of functions in $C_c(X)$ with compact support in $C$, the above calculation shows that $P_\gamma f=0$ except when this representation contains a $\gamma$-component. 
 \snl
 We have  that $P_\gamma f_{ij}=f_{ij}$ when $f_{ij}$ is defined as before by $f_{ij}(x)=\int f(x\cdot u)\gamma_{ij}(u\inv),du$. Indeed, the map $P_\gamma$ projects $\mathcal F$ onto the $\pi$-invariant subspace obtained as the component of $\pi$ that is equivalent with $\gamma$. 
\opm\label{opm:4.8}
We can not work with unitary representations as we have no  scalar product on the representation space $\mathcal F$. For this reason, the orthogonality properties of the matrix elements are expressed as
\begin{equation*}
\int \gamma_{ij}(u)\gamma_{k\ell}(u\inv)\,du=n_\gamma\delta_{i\ell}\delta_{jk}.
\end{equation*}
In the case of a unitary representation, when we use an orthonormal basis, we have $\gamma_{k\ell}(u\inv)=\overline \gamma_{\ell k}(u)$ and the formula reads as the more common version of the orthogonality of matrix elements.
\eopm

Here we do not work with unitary representations and this explains the appearance of $u\inv$ in these expressions.

 \prop\label{prop:4.9} 
 For all $\gamma$ we have $P_\gamma(\pi(u)f)=\pi(u)(P_\gamma f)$ for all $u\in G_0$ and $f\in\mathcal F$.
 \eprop
 \bew
 Let $f$ be in $\mathcal F$ and $\gamma$ any irreducible representation of the subgroup $G_0$. Then for all $v\in G_0$ we have
\begin{align*}
(P_\gamma (\pi(v) f)(x)
&=\int _{G_0} (\pi(v)f)(x\cdot u)\chi_\gamma(u\inv)\,du\\
&=\int_{G_0} f(x\cdot uv)\chi_\gamma(u\inv)\,du\\
&=\int_{G_0} f(x\cdot vv\inv uv)\chi_\gamma(u\inv)\,du\\
&=\int_{G_0} f(x\cdot vu)\chi_\gamma(v\inv u\inv v)\,du\\
&=\int_{G_0} f(x\cdot vu)\chi_\gamma(u\inv)\,du\\
&=(\pi(v)P_\gamma f)(x)
\end{align*}
So $\pi(v)$ and $P_\gamma$ commute for all $v\in G_0$.
\ebew

\nl
\bf Polynomial functions and different compact open subgroups \rm
\nl
Suppose that we are given two polynomial functions $f_1$ and $f_2$. Then we have compact subsets $C_1$ and $C_2$, containing the support of $f_1$ and $f_2$ resp.\ and two compact open subgroups $G_1$ and $G_2$ with the property that $C_1\cdot G_1=C_1$ and $C_2\cdot G_2=C_2$. The function $f_1$ belongs to a finite-dimensional $G_1$-invariant subspace of $\mathcal F_{C_1}$ and $f_2$ to a finite-dimensional $G_2$-invariant subspace of $\mathcal F_{C_2}$.
\ssnl
Define $C=C_1\cup C_2$ and $G_0=G_1\cap G_2$. The set $C$ is again a compact subset of $X$, also $G_0$ is still a compact open subgroup of $G$ and we have that $C\cdot G_0=C$. Now the functions $f_1$ and $f_2$ both belong to a finite-dimensional $G_0$-invariant subspace of functions in $\mathcal F_C$. See item ii) in  Remark \ref{opm:3.17}.
\ssnl
We now concentrate on one of the inclusions. More precisely we consider a compact subset $C$ of $X$, a compact open  subgroup  $G_0$ of $G$ satisfying $C\cdot G_0=C$ and another compact open subgroup $G_1$ of $G$, contained in $G_0$. Also for this subgroup $G_1$ we have $C\cdot G_1=C$.
\ssnl
A finite-dimensional $G_0$-invariant subspace of functions with compact support in $C$ will also be $G_1$-invariant. In particular, given $\gamma$ and $u\in G_1$ we have $P_\gamma (\pi(u)f)=\pi(u) (P_\gamma f)$. Finite-dimensional subspaces of $\mathcal F$, invariant under the action of $G_0$ will split into $G_1$-invariant subspaces.  
\ssnl
More precisely we have the following.
\prop\label{prop:4.10}
Assume that $f$ is a function with compact support in $C$ and that we can write
\begin{equation*}
f(x\cdot u)=\sum_i f_i(x)\varphi_i(u)
\end{equation*}
for all $x$ and $u\in G_1$, 
where the $f_i$ again have support in $C$ and the $\varphi_i\in C_c(G)$. Then we also have
\begin{equation*}
(P_\gamma f)(x\cdot u)=\sum_i (P_\gamma f_i)(x)\varphi_i(u)
\end{equation*}
for all $x$ and all $u\in G_1$.
\eprop

\bew
Suppose that $f$, with support in $C$, is polynomial for the action of $G_0$. Write $f(x\cdot u)=\sum_i f_i(x)\varphi_i(u)$. Then
\begin{equation*}
(P_\gamma f)(x\cdot u)=P_\gamma \pi(u) f=\sum_i P_\gamma f_i (x)\varphi_i(u).
\end{equation*}
We have used Proposition \ref{prop:4.9}.
\ebew

\nl
\bf Existence of enough polynomial functions \rm
\nl
In the first place, we have the following result.

\prop\label{prop:4.11}
If $C$ is  a compact subset of $X$ such that $C\cdot G_0=C$, then the polynomial functions in $C_c(X)$ with  support in $C$ are dense in $\mathcal F_C$ for the supremum norm.
\eprop

\bew
Given $f$ in $C_c(X)$ with compact support in $C$ we have the functions $f_{ij}$ as defined in Notation \ref{notat:4.4}, for all irreducible representations $\gamma$ of $G_0$. From the representation theory of compact groups, we know that the matrix elements of irreducible representations span the subalgebra $\mathcal P(G_0)$ and it is dense in the function algebra $C(G_0)$. Then we can approximate the function $f$ by polynomial functions of the form 
\begin{equation*}
x\mapsto \int_{G_0} f(x\cdot u)\varphi(u\inv)\,du
\end{equation*}
where $\varphi\in \mathcal P(G_0)$.  
\ebew

We have at least two special cases where Proposition \ref{prop:4.11} can be used to provide polynomial functions. 
\ssnl
First there is the case of a discrete group $G$. Then the trivial subgroup $\{e\}$, only containing the identity element $e$, 
is a compact open subgroup and all compact subsets $C$ of $X$ can be used. This however does not give anything new because for a partial action of a discrete group, all functions $f\in C_c(X)$ are polynomial. See Remark \ref{opm:3.4}. 
\ssnl
On the other hand, when we have a \emph{global} action of a locally compact group with a compact open subgroup on a locally compact space, the above result shows that the polynomial functions are dense in $C_c(X)$.
\ssnl
For a global action we can also prove the following. 

\prop\label{prop:4.12}
Assume that we have a global action of a group $G$ with a compact open subgroup $G_0$. 
Given a compact set $C$ in $X$ there exist a polynomial function $g$ on $X$ so that $g(x)=1$ for all $x\in C$.

\eprop
\bew
As in the proof of the Proposition \ref{prop:4.10} above, start with any function $f\in C_c(X)$ and a polynomial function $\varphi$ on $G_0$ and define $g$ again as
\begin{equation*}
g(x)=\int f(x\cdot u)\varphi(u)\,du.
\end{equation*}
We know that $g$ is a polynomial function on $X$. Now suppose that $f(x\cdot h)=1$ for all $x\in C$ and $h\in G_0$. Such a function exist because $C$ and $G_0$ are compact.
 Then, for $x\in C$ we have
\begin{equation*}
g(x)=\int f(x\cdot u)\varphi(u\inv)\, dh=\int \varphi(u\inv)\,du.
\end{equation*}
If we normalize $\varphi$ so that $\int \varphi(u)\,du=1$ we get $g(x)=1$ for $x\in C$.
\ebew

In fact, we could have taken for $\varphi$ the constant function $1$ on $G_0$. And we see that we will be able to find $g$ so that not only $g(x)=1$ on $C$ but also that $0\leq g\leq 1$. To achieve this, we just need to assume that this last property holds also for $f$.
\ssnl
It is a consequence of this result that when $G$ has a compact open subgroup,  the algebra $C_c(X)$ with pointwise multiplication, has local units in the subalgebra $\mathcal P(X)$ of polynomial functions when $G$. 

\opm\label{opm:4.13}
We see from Proposition \ref{prop:4.10} that, even when the group $G$ has a compact open subgroup, the existence of enough polynomial functions is only guaranteed when there are enough compact subsets that are invariant under the action of a compact open subgroup. This is a fairly strong condition but as we will see in the next sections, there are some interesting cases where this happens. We also refer to \cite{La-VD3} where fortunately we  have such actions.
\eopm

We have made such a remark already earlier in the paper, see after Corollary  \ref{gev:3.19a}.
\ssnl
Now we give an example to illustrate this remark. 

\voorb\label{voorb:4.15} 
i) First consider the global right action $(x,p,q)\mapsto (x,pq)$ of $G$ on $X\times G$  as in Proposition \ref{prop:1.15}. We have non-trivial polynomial functions if and only if  $G$ has a compact open subgroup. 
\ssnl
ii) Next restrict this global action to the open set $\Gamma$. Then we are interested in  polynomial functions with compact support in $\Gamma$. Even if the group $G$ has a compact open subgroup, it is not obvious that such polynomial functions exist. 
\ssnl
iii) If we do have such functions, we will also have polynomial functions on $\Gamma$ for the action $(x,p)\tl_1 q= (x\cdot q,q\inv p)$, as  considered in Proposition \ref{prop:1.14}. This is so because the map $(x,p)\mapsto (x\cdot p,p\inv)$ intertwines these actions. In the end, this will provide polynomial functions for the original action of $G$.
\ssnl
iv)
Also conversely, if we have non-trivial polynomial functions for the original partial right action, we are able to construct polynomial functions on $\Gamma$ for the partial action of $G$ and eventually, we get polynomial functions on $X\times G$ for the global right action, but with compact support in the subset $\Gamma$. 
\evoorb

Related properties for this case are found in Proposition \ref{prop:6.6} 
in Section \ref{s:misc}.

\section{\hspace{-17pt}. Polynomial functions on a locally compact group}\label{s:pol-lcg}  

In this section, we let a locally compact group $G$ act on itself by multiplication from the right. We study polynomial function on $G$ for this action. We recover results also found in \cite{La-VD1} from a different point of view. Some new results are included. We use the definition of a polynomial function as given in this paper, see Definition \ref{defin:3.1}.
\ssnl
We not only consider right multiplication but also the action given by left multiplication and the combination of the two. An important issue is the existence of polynomial functions with compact support in a given open subset of $G$. At the end in this section, we have an item on local units for the function algebra $C_c(G)$ as well as for the convolution algebra $C_c(G)$.
\ssnl
We have seen in a previous section that polynomial functions on $G$ for right multiplication can only exist when the group $G$ has a compact open subgroup, see Corollary \ref{gev:3.19a}. Therefore, it makes sense to start with the easier and almost trivial case of polynomial functions on a compact group.
\nl
\bf Polynomial functions on a compact group \rm
\nl
Assume that $G$ is a compact group. We let it act on itself by right multiplication. This is a global action. It is proper and locally homeomorphic as we already mentioned in  remarks after Definition \ref{defin:1.5} and Definition \ref{defin:1.6}.
\ssnl
We consider the representation $\pi$ of $G$ on $C(G)$ given by $(\pi(q)f)(p)=f(pq)$. Then we have the following characterization of polynomial functions on $G$. 

\prop\label{prop:5.1}
A function $f$ is polynomial if and only if it belongs to some finite-dimensional $\pi$-invariant subspace of $C(G)$.
\eprop

\bew
i) Let $F$ be any finite-dimensional invariant subspace of $C(G)$ and let $f\in F$. Then, because $F$ is an invariant subspace, we have $p\mapsto f(pq)$ in $F$ for all $q$. We can take for $V$ in Definition \ref{defin:3.1} the whole group $G$. This proves that $f$ is polynomial.
\ssnl
ii) Conversely let $f$ be polynomial on $G$. Choose a neighborhood $V$ of the identity and write
$f(pv)=\sum_i\varphi_i(v)f_i(p)$ for $v\in V$ as in Proposition \ref{prop:3.8}. 
We can cover $G$ by finitely many sets $\{q_kV\}_k$ and take a partition of unity $\{\lambda_k\}_k$ in $C(G)$ for this cover. Then we have for all $p,q\in G$ that
\begin{align*}
f(pq)&=\sum_k \lambda_k(q)f(pq)\\
&=\sum_k \lambda_k(q)f(pq_kq_k\inv q)\\
&=\sum_{ik} \lambda_k(q)\varphi_i(q_k\inv q)f_i(pq_k).
\end{align*}
We have used that $q_k\inv q\in V$ if $\lambda_k(q)\neq 0$.
Redefining the functions $f_i$, we see that we can write $f(pq)=\sum f_i(p)g_i(q)$, now for all $p,q\in G$. Assume that the functions $g_i$ are linearly independent and then let $F$ be the subspace of $C(G)$ spanned by  the functions $f_i$. With $q=e$ we see that $f\in F$. Further
\begin{align*}
\sum_i f_i(pr)g_i(q)=f(prq)=\sum_i f_i(p)g_i(rq)
\end{align*}
for all $p,q,r$ and because the functions $g_i$ are supposed to be linearly independent, we can conclude that the functions $p\mapsto f_i(pr)$ again belong to $F$. We can use Proposition \ref{prop:3.6} as in the proofs of Propositions \ref{prop:3.7} and \ref{prop:3.8} to make this precise. 
\oldcomment{Check internal references! This is done. It may look as an unnecessary remark, but it is not.}
\ssnl
Hence $F$  is a  finite-dimensional $\pi$-invariant subspace of $C(G)$ containing $f$.
\ebew

An argument based on Proposition \ref{prop:3.18} would give us a compact open subgroup $G_0$ of $G$. The technique used  in item ii) in the proof tells us how to lift the result to all of $G$. In fact, the same idea is used in Proposition \ref{prop:6.1} and further in Proposition \ref{prop:6.4a} and Proposition \ref{prop:6.14a}.
\ssnl
 With the notations from the proof above, we have $f=\sum g_i(e)f_i$ and when we write
\begin{equation}
f_i(pq)=\sum_j \pi_{ji}(q)f_j(p)\label{eqn:4.2}
\end{equation}
we get $f_i(q)=\sum_jf_j(e)\pi_{ji}(q)$. So $f(q)=\sum_{ij} g_i(e)f_j(e)\pi_{ji}(q)$. We see that $f$ is a linear combination of matrix elements in the subrepresentation of $G$ on $F$.
\ssnl
In fact, it is easy to see that such matrix elements (as functions on the group) are always polynomial.
\ssnl
From the representation theory of compact groups, we know that the representation $\pi$ of $G$ on $C(G)$ is a direct sum of finite-dimensional representation. It follows that there are plenty of polynomial functions in the sense that they are dense in $C(G)$. See e.g.\ Chapter 5 in \cite{Fo}. 
\nl
The situation is completely different for \emph{non-compact groups}. Matrix elements in continuous finite-dimensional representations  will no longer belong to $C_c(G)$. So they are not polynomial in the sense of Definition \ref{defin:3.1}. We treat this case in the next subsection.
\nl
\bf Polynomial functions on a group with a compact open subgroup \rm
\nl
Let $G$ be a locally compact group. Also here we consider the global action of $G$ on itself given by right multiplication. 
As in the case of a compact group, we consider the representation $\pi$ of $G$ on $C_c(G)$ given by $(\pi(q)f)(p)=f(pq)$. We now get the following characterization of polynomial functions.

\prop\label{prop:5.2}
A non-zero function $f$ in $C_c(G)$ is polynomial if and only if there is a compact open subgroup $G_0$ of $G$ so that $f$ belongs to a finite-dimensional $G_0$-invariant subspace.
\eprop

\bew
i) Suppose that $G_0$ is a compact open subgroup and $F$ a finite-dimensional  $G_0$-invariant subspace. For $f$ in $F$ we have all the functions $p\mapsto f(ph)$ in $F$ for all $h\in G_0$. Now we can use $G_0$ for the neighborhood of the identity needed in Definition \ref{defin:3.1}.This proves that $f$ is polynomial.
\ssnl
ii) Conversely let $f$ be polynomial in $C_c(G)$. We know from Theorem \ref{stel:3.15} and Proposition \ref{prop:3.18} that there is a compact open subgroup $G_0$ of $G$ and a finite-dimensional subspace $F$ of $C_c(G)$, containing $f$ and invariant under the action of elements in $G_0$. This precisely means that $F$ is a finite-dimensional $G_0$-invariant subspace. This proves the result. 
\ebew

We see that the two characterizations of polynomial functions, in Proposition \ref{prop:5.1} for the compact case and Proposition \ref{prop:5.2} for the general locally compact case, are essentially straightforward  reformulations of general results obtained in  Section \ref{s:pol}. The technique used in item ii) of the proof of Proposition \ref{prop:5.1} can be used to show that the result above will hold for any compact open subgroup of $G$. In particular, we have the following.

\prop\label{prop:5.2a}  
A function $f\in C_c(G)$ is polynomial in the sense \cite{La-VD1} if and only if it is polynomial in the sense of Definition \ref{defin:3.1}. 
\eprop

\bew
First recall the definition from \cite{La-VD1}. We refer to Definition 3.6 and Lemma 3.3 of that paper. A function $f$ is polynomial if and only if there is a compact open subgroup $G_0$ of $G$ and finitely many funictions $f_i\in C_c (G)$ and $\varphi_i\in C(G_0)$ satisfying $f (pv) = \sum_i f_i(p)\varphi_i(v)$ for all $p\in G$ and $v\in G_0$.
\ssnl
i) If $f$ is polynomial in the sense of \cite{La-VD1}, we can use $G_0$ for the neighborhood $V$ in our Definition \ref{defin:3.1}. Recall that for this global action, we have $V_f=G$ for all $f$.
\ssnl
ii) Conversely, if $f$ is polynomial in the sense of our Definition \ref{defin:3.1}, we can use Proposition \ref{prop:3.13} in combination with Proposition \ref{prop:3.18}.  
\ebew

We will come back to this result in Section \ref{s:misc}, after the proof of Proposition \ref{prop:6.2}. \keepcomment{We may eventually omit this proof}.
\nl
Now  from the representation theory of compact groups (see e.g\ \cite{Fo}), we know that the restriction of the representation of $G$ on $C_c(G)$ to a compact open subgroup $G_0$ will be a direct sum of finite-dimensional representations. This again implies the existence of plenty of polynomial functions. It will imply that the space of polynomial functions in $C_c(G)$ is also dense with respect to the sup norm when $G$ has a compact open subgroup.
\ssnl
In fact, Proposition \ref{prop:4.11} for functions on $G$ provides an obvious way to construct polynomial functions on $G$ using polynomial functions on $G_0$. 

\prop\label{prop:5.3}
Let $f$ be any function in $C_c(G)$. Assume that $G_0$ is  a compact open subgroup of $G$ and that $\varphi$ is a polynomial function on $G_0$. Define $g$ on $G$ by
\begin{equation*}
g(p)=\int_{G_0} f(ph)\varphi(h)\,dh
\end{equation*}
where we integrate over the Haar measure on $G_0$. Then $g$ is a polynomial function on $G$.
\eprop

We also have the following property. We include a proof for completeness.

\prop\label{prop:5.4}
Assume that the group $G$ has  a compact open subgroup $G_0$. 
Given a compact set $C$ in $G$ there exist a polynomial function $g\in C_c(G)$ so that $g(x)=1$ for all $x\in C$.
\eprop

\bew
As in the proof of the Proposition \ref{prop:5.3} above, start with any function $f\in C_c(G)$ and a polynomial function $\varphi$ on $G_0$ and define $g$ as
\begin{equation*}
g(x)=\int f(ph)\varphi(h)\,dh
\end{equation*}
for $p\in G$. 
We know that $g$ is a polynomial function on $G$. Now suppose that $f(p h)=1$ for all $p\in C$ and $h\in H$. Such a function exist because $C$ and $G_0$ are compact.
 Then, for $p\in C$ we have
\begin{equation*}
g(x)=\int f(ph)\varphi(h)\, dh=\int \varphi(h)\,dh.
\end{equation*}
If we normalize $\varphi$ so that $\int \varphi(p)\,dp=1$ we get $g(p)=1$ for $p\in C$.
\ebew 

In fact, we could have taken for $\varphi$ the constant function $1$ on $G_0$. And we see e.g.\ that we will be able to find $g$ so that not only $g(x)=1$ on $C$ but also that $0\leq g\leq 1$. To achieve this, we just need to assume that this last property holds also for $f$.
 
 \opm
 Let $G_0$ be a compact open subgroup of $G$.  Assume that $C$ is a compact subset of $G$. 
 We can cover it  by finitely many translates $(p_kG_0)_{k=1}^n$, where $p_k\in C$ for all $k$. Then we can define a map $T$ from $C_c(G)$ to the  direct sum of $n$ copies of $C(G_0)$ by $(T(f))_k(h)=f(p_kh)$. If $T(f)=0$ we must have $f(p)=0$ for all $p\in C$. The map $T$ also intertwines the representation of $G_0$ on $C_c(G)$ with the direct sum of $n$ copies of the right regular representation of $G_0$ on $C(G_0)$. Furthermore we know that any finite-dimensional representation of $G_0$ can occur only with a finite multiplicity in the regular representation. It follows from this that the same is true for the space of functions with a given compact support.
 \eopm
 
 The following result essentially says the same thing.
 
 \prop\label{prop:5.6}
 Consider a compact subset $C$ of $G$. 
 Assume that we are given finitely many functions $(\varphi_i)$ in $C(G_0)$. Consider the space $F$ of  functions $f$ with compact support in $C$ so that there exist functions $f_i$  satisfying 
 \begin{equation*}
f(ph)=\sum_j \varphi_j(h)f_i(p)
\end{equation*}
for all $p\in C$ and $h\in G_0$. Then $F$ is finite-dimensional.
 \eprop
 
 \bew
 Cover $C$ with finitely many sets $p_kG_0$ where $p_k\in C$.  Let $(\gamma_k)$ be a partition of unity in $C_c(G)$ for this covering. Then we can write
 \begin{equation*}
f(p)=\sum_k \gamma_k(p)f(p)=\sum_k \gamma_k(p)f(p_kp_k\inv p)
\end{equation*}
for all $p$. By the presence of $\gamma_k(p)$, the element $p_k\inv p$ is forced to lie in $G_0$ and therefore we find
\begin{equation*}
f(p)=\sum_{ik}\gamma_k(p) \varphi_i(p_kp)f(p_k).
\end{equation*}
We see that $f$ belongs to the space  spanned by the functions $p\mapsto \gamma_k(p)\varphi_i(p_kp)$. This proves the result.
 \ebew
 
 A similar argument is used in the proof of Proposition \ref{prop:5.1} and again further in Proposition \ref{prop:6.14a}.
 
\nl
\bf Further properties of polynomial functions on a locally compact group \rm
\nl
\keepcomment{We have to check the logical structure of this part here with references to the first part of Section 6.
\ssnl
Also refer to the multiplier Hopf algebras sitting inside $C_c(G)$. We do this in Section 6.
\nl}{}
We now prove a few extra properties for elements in $\mathcal P(G)$, also found in the original paper \cite{La-VD1}, but treated here slightly differently.

\prop 
Suppose that $f\in \mathcal P(G)$. Let $q\in G$, then the function $p\mapsto f(qp)$  is again polynomial. 
\eprop

This is  an immediate consequence of Definition \ref{defin:3.1}. We can use the same neighborhood $V$ as for $f$ itself. Then, for all $v\in V$, the function $p\mapsto f(qpv)$ will belong to the space of left translations of the functions in the original finite-dimensional space.
\ssnl
In the case of a group, we can also consider left multiplication and develop polynomial functions for this action. From the following result, it will follow that this gives the same class of functions.

\prop
Assume that  $f$ is polynomial function on $G$. Then for all compact subsets $C$ of $G$, there is a finite-dimensional subspace $F$ of $C_c(G)$ containing all the functions $p\mapsto f(qp)$ with $q\in C$.
\eprop

\bew
Denote the function $p\mapsto f(qp)$ by $f_q$. Because $f$ has compact support and $C$ is compact, we can find a compact subset $C_1$ of $G$ so that all the functions $f_q$, for $q\in C$, have support in $C_1$. Because $f$ is polynomial we have functions $\varphi$ and $f_i$ so that 
\begin{equation*}
f_q(ph)=f(qph)=\sum_i \varphi_i(h)f_i(qp).
\end{equation*}
This holds for all $q\in C$ and all $p\in G$. Then we can apply Proposition \ref{prop:5.6} and we find that there is a finite-dimensional subspace $F$ of $C_c(G)$ containing all the functions $f_q$ with $q\in C$.
\ebew

We can reformulate this result in the following way.

\gev\label{gev:5.9}
Given a polynomial function $f$ on $G$ and any $g\in C_c(G)$, we have functions $g_i$ and $f_i$ in $C_c(G)$ so that
\begin{equation*}
g(q)f(qp)=\sum_i g_i(q)f_i(p)
\end{equation*}
for all $p,q$. When the functions $g_i$ are chosen to be linearly independent, then the functions $f_i$ will again be polynomial.
\egev

In the next section, we will prove a generalization of this result, see Proposition \ref{prop:6.1}. In fact, also for a local action, we have a similar property, see Proposition \ref{prop:6.4a}.

\opm 
i) We can conclude that a polynomial function on $G$ for right multiplication, is also polynomial for  left multiplication.  And by symmetry, the converse is true as well.
\ssnl
ii) We also have this result on the other side. Given a polynomial function $f\in C_c(G)$ and any $g\in C_c(G)$ we can write $f(pq)f(q)=\sum_i f_i(p)g_i(q)$. 
\eopm

As an easy consequence, we get the following.

\prop\label{prop:5.11}
If $f$ is a polynomial function, then also $p\mapsto f(p\inv)$ is a polynomial function.
\eprop
\bew
Denote $p\mapsto f(p\inv)$ by $\widetilde f$. 
By the above result, given a compact subset $C$ of $G$, we have  a finite-dimensional subspace of $C_c(G)$ containing all the functions $p\mapsto f(qp)$ with $q\in C$. If we take a symmetric neighborhood  $V$ of $e$ with compact support, we see that the functions  $p\mapsto \widetilde f(pv)$ all belong to a finite subspace of $C_c(G)$. 

\ebew

It is interesting to compare the proof of this property with the one given in the original paper, see Theorem 3.7 in \cite{La-VD1}.
\ssnl

If we combine these results, we  obtain the following property.

\prop \label{prop:5.12}
Consider the right action of $G^\text{op}\times G$ on $G$ given by $p\cdot (q_1,q_2)=q_1p q_2$. If $f$ is  a polynomial function on $G$ (for the right action of $G$), then it is also polynomial for this action of $G^\text{op}\times G$. 
\eprop

\bew
Because $f$ is polynomial, we have a neighborhood $V$ of the identity and functions $f_i$ and $\varphi_i$ so that $f(pv)=\sum_i \varphi_i(v)f_i(p)$ for all $p\in G$ and $v\in V$. We also know from Proposition \ref{prop:3.8}
that we can  assume that the functions $f_i$ are again polynomial functions. As we have seen, they also are polynomial for left multiplication. Therefore, for each $i$ we can find a neighborhood $V_i$ of the identity and functions ${g_{ij}}$ and ${\psi_{ij}}$ so that $f_i(wp)=\sum_{i,j} \psi_{ij}(w)g_{ij}(p)$ for all $w\in V_i$ and all $p$. Now let $v\in V$ and $w$ in the the intersection of all the $V_i$. Then we find for all $p\in G$
\begin{equation*}
f(wpv)=\sum_i \varphi_i(v)f_i(wp)=\sum_{i,j}\varphi_i(v)\psi_{ij}(w)g_{ij}(p).
\end{equation*}
This proves that $f$ is polynomial for the action of $G^\text{op}\times G$.
\ebew

We have an example of commuting actions. There is right multiplication by $G$ on $G$ and left multiplication. The left action gives the right action of $G^\text{op}$ and then they combine into the right action of $G^\text{op}\times G$ on $G$ as in Proposition \ref{prop:1.3} of Section \ref{s:loc-act}.
\ssnl
We have here a special case of some other results we will have in the next section. First when $f$ is a polynomial function on $G$ for right multiplication, it is in the first place, also polynomial for left multiplication. Then the above result is a consequence of the more general result we obtain in Proposition \ref{prop:6.17a} in the next section. Secondly it would also be a consequence of the result of Proposition \ref{prop:6.15a} where we take $G^\text{op}\times G$ for the larger group and $G$ for the subgroup. We indeed have a situation where a function is polynomial as soon as it is already polynomial for the action of a smaller group.
 
 \nl
\bf Two algebras with underlying space $\mathcal P(G)$ \rm
\nl
The vector space $C_c(G)$ can be made into a $^*$-algebra in two ways. First one has the pointwise product and then the involution is simply given by $f\mapsto \overline f$ where $\overline f(p)=\overline{f(p)}$. We refer to this as the \emph{function algebra} $C_c(G)$. The other one is the \emph{convolution algebra} $C_c(G)$. The product is defined by
\begin{equation*}
(f*g)(p)=\int f(q)g(q\inv p)\,dq
\end{equation*}
where we integrate over the left Haar measure of $G$. The involution here is given by $f\mapsto f^*$ where $f^*(p)=\delta(p)\inv \overline{f(p\inv)}$. Here $\delta$ is the modular function of $G$.
\ssnl
It is obvious that the pointwise product of polynomial functions is again a polynomial function. In fact, we have the following result, essentially obtained already in Proposition \ref{prop:4.12}.

\prop 
Assume that $G$ has a compact open subgroup. 
The space $\mathcal P(G)$ of polynomial functions on $G$ is a non-degenerate $^*$-subalgebra of the function algebra $C_c(G)$. It separates points of $G$. For each compact subset $C$ of $G$ there is a function $f\in \mathcal P(G)$ so that $f(p)=1$ for all $p\in C$. Consequently the function algebra $C_c(G)$ has local units in the subalgebra $\mathcal P(G)$. Therefore  the product in the algebra $\mathcal P(G)$ is still non-degenerate.
\eprop

Remark that a $^*$-subalgebra of the function algebra $C_c(G)$ is always non-degenerate. Indeed, if $g$ is an element in the subalgebra satisfying $fg=0$ for all $f$ in the subalgebra, then $g^*g=0$ and hence $g=0$.
\ssnl

If there is no compact open subgroup, there are no non-zero polynomial functions and the result is no longer true. But of course, the algebra still has local units as such.
\nl

Next we look at $C_c(G)$ with the {\emph{convolution product}. 
\ssnl 
Because of the presence of the modular function in the formula for the involution here, we first prove the following lemma.

\lem\label{lem:5.14}
Assume that $G$ has a compact open subgroup. Let $\delta$ be any continuous homomorphism from $G$ to $\mathbb R^+$.
Then $\delta$ is locally constant on $G$. On a compact subset of $G$ it can take only finitely many values.
\elem

\bew
Let $G_0$ be a compact open subgroup. Then $\delta(h)=1$ for all $h\in G_0$ because $G_0$ is a compact subgroup of $G$ and $\delta$ is a continuous homomorphism from $G$ to $\mathbb R^+$.
For all $p$ we have a neighborhood $pG_0$ of $p$ and $\delta(ph)=\delta(p)$ because $\delta$ is a homomorphism. This proves that $\delta$ is locally constant.
\ssnl
Let $C$ be a compact subset of $G$. Then $C$ is contained in the union of finitely many
cosets $pG_0$ . And because $\delta$ equal to $1$  on $G_0$, we see that
the image of $C$ under $\delta$ is a finite subset of $\mathbb R^+$ . 
\ebew

We use this in the proof of the following proposition.

\prop Assume that $G$ has a compact open subgroup.
Let $f,f'$ be in $C_c(G)$ and denote by $f*f'$ the convolution product. Then $f*f'$ is polynomial as soon as one of the functions is polynomial. So $\mathcal P(G)$ is a two-sided $^*$-invariant ideal in the convolution algebra $C_c(G)$. 
It is a $^*$-subalgebra  with local units.  The product is still non-degenerate.
\eprop
\bew
i) First we show that the convolution of  functions is again polynomial as soon as one of them is polynomial. For this, take $f,f'\in C_c(G)$ and assume that $f'$ is polynomial. Let $V$ be a neighborhood of the identity in $G$ so that the functions $p\mapsto f'(pv)$ belong to a finite-dimensional space $F$. For all $p$ in $G$ and all $v$ in $V$ we have
\begin{equation*}
(f*f')(pv)=\int f(q)f'(q\inv pv)\,dq
\end{equation*}
If we let $F'$ be the space of functions $f*g$ where $g\in F$ we see that $(f*f')(\,\cdot\, v)$ belongs to $F'$ for all $v\in V$. So $f*f'$ is again polynomial. A similar argument will show that $f*f'$ is polynomial if $f$ is polynomial.
\ssnl
ii) Also $f^*$ is polynomial when this is the case for $f$. 
In Proposition \ref{prop:5.11} we have seen that $p\mapsto f(p\inv)$ is again polynomial. Taking the complex conjugate gives no problem. In Lemma \ref{lem:5.14} we have seen that $\delta$ can only take finitely many values on a compact set. 
Therefore  the map $p\mapsto \delta(p)\inv f(p\inv)$ will still be a polynomial function on $G$. 
\ssnl
iii) 
We now prove the existence of local units. 
Let $f$ be a polynomial function. Assume that $G_0$ is a compact open subgroup. We can write
\begin{equation*}
f(h\inv p)=\sum \varphi_j(h)f_j(p).
\end{equation*}
We can assume that the functions $f_j$ are linearly independent so that the functions $\varphi_j$ are again polynomial. Furthermore, we may assume that the functions $\varphi_j$ are mutually orthogonal in $C(H)$ and normalized so that $\int |\varphi_j(h)|^2\,dh=1$ for all $j$.  Put $g=\sum_j  \varphi_j(e)\overline{\varphi_j}$. It will be polynomial. Further we have 
\begin{align*}
(g*f)(p)
&=\int g(h)f(h\inv p),dh\\
&=\int \sum g(h)\varphi_j(h)f_j(p)\,dh\\
&=\sum \varphi_j(e)f_j(p)=f(p)
\end{align*}
 and we see that $g*f=f$. 
\ssnl
iv) In a similar way we find for all $f\in \mathcal P(G)$ an element $g\in \mathcal P(G)$ so that $f*g=f$. 
This guarantees the existence of  local units by a simple argument as found in \cite{Ve}. 

\ebew

Also here there is a more direct way to show that the subalgebra $\mathcal P(G)$ of the convolution algebra is still non-degenerate. Indeed, let $g$ be in $\mathcal P(G)$ and assume that $f*g=0$ for all $f$ in the $\mathcal P(G)$. Then $g^**g=0$. This implies $g=0$ because the convolution algebra is an operator algebra. 
\ssnl
Remark that in general the convolution algebra $C_c(G)$ itself does not have local units, only approximate units.
In fact, we have the following property of the convolution algebra $C_c(G)$.

\prop\label{prop:5.17}
Assume that $f,\xi$ are functions in $C_c(G)$ satisfying $f*\xi=\xi$. Then $\xi$ is a polynomial function on $G$. With $\xi\neq 0$, it follows that $G$ has a compact open subgroup.
\eprop
\bew
We have $\int f(q)\xi(q\inv p)\,dq=\xi(p)$ for all $p\in G$. Let $V$ be a neighborhood of the identity in $G$ with compact closure. Let $g$ be a function in $C_c(G)$ so that $g(pv)=1$ if $f(p)\neq 0$ and $v\in V$.  Denote by $\lambda(f)$  convolution with $f$  and by $\pi(g)$  multiplication with the function $g$. Define $\xi_v(p)=\xi(pv)$ for all $p$ and $v\in V$.
Then
\begin{align*}
(\pi(g)\lambda(f)\xi_v)(p)
&=g(p)\int f(q)\xi_v(q\inv p)\,dq\\
&=g(p)\int f(q)\xi(q\inv pv)\,dq\\
&=g(p)\xi(pv)=\xi(pv)=\xi_v(p).
\end{align*}
Hence we have $\pi(g)\lambda(f)\xi_v=\xi_v$ for all $v\in V$. Now we know that the operator $\pi(f)\lambda(g)$ is a compact operator. The eigenspaces are finite-dimensional and so there is a finite-dimensional subspace $F$ of $C_c(G)$ containing  $\xi_v$ for all $v\in V$. Therefore $\xi$ is polynomial.
\ebew

Here is another proof of this result. Denote by $T$ the operator $\pi(f)\lambda(g)$ on $C_c(G)$. It is given by
\begin{equation*}
(T\xi)(p)=\int \varphi(p,q)\xi(q\inv p)\,dp=\int\varphi(qp,q)\xi(p)\,dp
\end{equation*}
where $\varphi(p,q)=f(p)g(q)$. Let $(\xi_j)$ be an orthonormal basis for the space $F$ spanned by the functions $\xi_v$. Then we have
\begin{align*}
\sum_j \|\xi_j\|^2 = \sum_j \|T\xi_j\|^2
&=\sum_j\int |(T\xi_j)(p)|^2\,dp \\
&=\int \sum_j \left| \int \varphi(qp,q)\xi_j(q)\,dq\right|^2\, dp\\
&\leq \int \int |\varphi(qp,q)|^2\,dqdp.
\end{align*}
Because $\varphi$ has compact support, this integral is finite. Therefore there can only be finitely many elements $\xi_j$ and the space $F$ is finite-dimensional.
\nl
\bf Multiplier Hopf algebras inside $C_c(G)$ \rm
\nl
It has been shown in \cite{La-VD2} that the function algebra $C_c(G)$ can only contain a non-trivial multiplier Hopf algebra if $G$ has a compact open subgroup. The same is true for the convolution algebra. We review these results here based on Proposition \ref{prop:5.17}. The argument uses the existence of local units in a multiplier Hopf algebra. This property is found in the literature on multiplier Hopf algebras, but we include here a proof of this result. We give more comments and references in a remakr after the proof.

\prop\label{prop:5.18}
Let $(A,\Delta)$ be a multiplier Hopf algebra. Then $A$ has local units.
\eprop 
\bew
Take any $a\in A$ and assume that $\omega$ is a linear functional on $A$ so that $\omega(ba)=0$ for all $b\in A$. Then
\begin{equation*}
(\iota\ot\omega)((1\ot c)\Delta(d)(1\ot a))=0
\end{equation*}
for all $c,d$. Hence also
\begin{equation*}
\sum_{(d)}\omega(S(d_{(1)})d_{(2)}a)=0
\end{equation*}
for all $d$. Here $S$ is the antipode. We use that $\Delta(d)(1\ot a)$ belongs to $A\ot A$. Then $\varepsilon(d)\omega(a)=0$ for all $d$ and $\omega(a)=0$. Now $\varepsilon$ is the counit. It follows that $a\in Aa$. Therefore there is an element $e$ in $A$ so that $ea=a$.
\ssnl
Similarly we have for all $a\in A$ an element $e\in A$ satisfying $ae=a$.
\ssnl
It now follows from an argument in \cite{Ve} that $A$ has local units.
\ebew

Here are some comments on the history of this result.

\opm
The first result of this type is found in \cite{Dr-VD-Z}. In Proposition 2.2 of that paper, it is shown that any \emph{regular} multiplier Hopf algebra has left and right local units. In Proposition 2.6 of the same paper, it is proven that any regular multiplier Hopf algebra with integrals has two-sided local units. When that paper was published (1999), it was still open if these properties would still hold for any multiplier Hopf algebra.
\ssnl
In \cite{VD-Z},  it is shown in Proposition 1.2 of that paper that in any multiplier Hopf algebra $(A,\Delta)$, given an element $a\in A$, there are elements $e,f$ in $A$ satisfying $ea=af=a$. The proof is as in the first part of the proof above of Proposition \ref{prop:5.18}. This paper was written later than \cite{Dr-VD-Z} but appeared also in 1999. At that time, the result of \cite{Ve} was not yet available and so it was still open wheter or not any multiplier Hopf algebra had local units.
\ssnl
Later, the existence of local units has also been obtained for algebraic quantum hypergroups and weak multiplier Hopf algebras. In the first case, the property is found in Proposition 1.6 of \cite{De-VD} while for the second case, it is proven in Proposition 2.14 of \cite{VD-W}.
\eopm

Then we get the  following as a consequence of the result of Proposition \ref{prop:5.18}.
\prop
Assume that $A$ is a (non-trivial) multiplier Hopf algebra sitting in the convolution algebra $C_c(G)$. Then $G$ has a compact open subgroup and all elements in $A$ belong to $\mathcal P(G)$.
\eprop

\bew
Given $f\in A$ we have an element $g \in A $ satisfying  $g*f=f$. Then $f$ is a polynomial function.
\ebew

We see that in fact, we do not need that $A$ has local units, the first property in the proof of Proposition \ref{prop:5.18} would be enough.
\ssnl
The result is true also for the function algebra.
\prop
Assume that $A$ is a (non-trivial) multiplier Hopf algebra sitting in the function algebra $C_c(G)$. Then $G$ has a compact open subgroup and all elements in $A$ are polynomial functions.
\eprop

\bew
The dual of $A$ is a multiplier Hopf algebra sitting in the convolution algebra and we can apply the previous result. The underlying spaces of $A$ and $\widehat A$ coincide. Therefore also the functions in $A$ are polynomial.
\ebew

We finish this section with a remark.

\opm
i) If $A$ is a multiplier Hopf algebra sitting in $C_0(G)$, it must sit in $C_c(G)$ because it has local units. Then the dual sits in the convolution algebras and because that also has local units, we must have that it belongs to $\mathcal P(G)$. So also $A$ belongs to $\mathcal P(G)$.
\ssnl
ii) What if $A$ is a multiplier Hopf algebra sitting in $C_r^*(G)$. If the integral is finite on $A$ we can construct the dual in $C_0(G)$ and apply the previous result. But is this always the case? If $A$ is a multiplier Hopf $^*$-algebra sitting in $C_r^*(G)$, we can apply the results of \cite{DC-VD}, Theorem 2.10 and Corollary 2.11. The result is also found in the paper \cite{La-VD1}, see Theorem 5.6. 
\eopm

\section{\hspace{-17pt}. Miscellaneous results on polynomial functions} \label{s:misc} 

In this section we collect some more properties of polynomial functions and show their use for constructing polynomial functions in special situations.

\nl
\bf Global behavior of polynomial functions \rm
\nl

Polynomial functions are defined using a local property of the action. In our first topic in this section, we derive some global properties. We prove the following result first for a global action. 

\prop\label{prop:6.1}
Suppose that we have a \emph{global} right action of $G$ on $X$ and a polynomial function $f$ in $C_c(X)$. Then for every $g\in C_c(G)$ we have finitely many functions $f_i\in C_c(X)$ and $g_i\in C_c(G)$ satisfying
\begin{equation}
f(x\cdot p)g(p)=\sum_i f_i(x)g_i(p)\label{eqn:6.1}
\end{equation}
for all $x\in X$ and $p\in G$. The functions $f_i$ can be chosen to be polynomial. Conversely, if we have a function $f\in C_c(X)$ such that for some non-zero $g\in C_c(G)$ there exists functions $f_i\in C_c(X)$ and $g_i\in C_c(G)$ satisfying Equation (\ref{eqn:6.1} above, then $f$ is polynomial.
\eprop

\bew
i) First assume that $f$ is polynomial. As in Proposition \ref{prop:3.8}, let $V$ be a neighborhood of the identity and let   $\varphi_i\in C_c(G)$ and $f_i\in C_c(X)$  satisfying
\begin{equation*}
f(x\cdot v)=\sum_i \varphi_i(v)f_i(x)
\end{equation*}
for all $x\in X$ and all $v\in V$. We can assume that the functions $f_i$ are polynomial (see Proposition \ref{prop:3.8}). 
\ssnl
Now take $g\in C_c(G)$ and let $C$ be its support. We can choose  elements $q_1, q_2, \dots, q_n$ in $C$ so that
\begin{equation*}
C\subseteq q_1 V \cup q_2 V \cup \dots \cup q_n V.
\end{equation*}
Choose a partition of unity $\lambda_1,\lambda_2,\dots\lambda_n$ of functions in $C_c(G)$, relative to this covering of $C$. We have for all $q$ in $G$ that
\begin{align*}
f(x\cdot q)g(q)
&=\sum_k \lambda_k(q)f(x\cdot q)g(q)\\
&=\sum_k \lambda_k(q) f(x\cdot q_k \cdot (q_k\inv q))g(q)\\
&=\sum_{i,k}  \lambda_k(q) f_i(x\cdot q_k)\varphi_i(q_k\inv q)g(q).
\end{align*}
For the last equality we use that $q_k\inv q\in V$ when $\lambda_k(q)\neq 0$. This means that we have
\begin{equation*}
f(x\cdot q)g(q)=\sum_{i,k}f_i(x\cdot q_k)\psi_{ik}(q)
\end{equation*}
where $\psi_{ik}(q)= \lambda_k(q) \varphi_i(q_k\inv q)g(q)$ for all $q\in G$. This proves the first statement. Remark that the functions $x\mapsto f_i(x\cdot q_k)$ are all still polynomial, see Proposition \ref{prop:3.7}.
\ssnl
ii) Conversely, take $f\in C_c(X)$, $g\in C_c(G)$ and assume that $g$ is non-zero. Take a point $q_0$ in $G$ with $g(q_0)\neq 0$ and a neighborhood $V$ of the identity in $G$ so that still $g(q)\neq 0$ for $q\in Vq_0$. From Equation (\ref{eqn:6.1}) 
we get
\begin{equation*}
f(x\cdot vq_0)g(vq_0)=\sum_i f_i(x)g_i(vq_0).
\end{equation*}
We divide by $g(vq_0)$ and we get
\begin{equation*}
f(x\cdot vq_0)=\sum_i f_i(x)g_i(vq_0)g(vq_0)\inv
\end{equation*}
when $v\in V$. This proves that $x\mapsto f(x\cdot q_0)$ is polynomial. Consequently, again by Proposition \ref{prop:3.7}, $f$ is polynomial as a translate of this function.
\ebew

We can apply this result for polynomial functions on $G$ and find the following. 

\prop\label{prop:6.2}
Let $G$ act on itself by right multiplication. If $f$ is a polynomial function on $G$ and $g$ any function in $C_c(G)$, we have functions $f_i$ and $g_i$ in $C_c(G)$ satisfying 
$$f(pq)g(q)=\sum_i f_i(p)g_i(q)$$
 for all $p,q\in G$. The functions $f_i$ can be chosen to be polynomial. Conversely, if $f\in C_c(G)$ and if $g$ is a non-zero function in $C_c(G)$ such that 
 $$f(pq)g(q)=\sum_i f_i(p)g_i(q)$$
  for all $p,q\in G$, then $f$ is polynomial.
\eprop

We recover the result from Corollary \ref{gev:5.9} as a special case of Proposition \ref{prop:6.1}.
\ssnl
This result implies that $C_c(G)$ can only contain a multiplier Hopf algebra if there are polynomial functions and therefore only when $G$ has a compact open subgroup. Hence we recover one of the main results proven in \cite{La-VD1}.
\nl
Now we consider the result for a \emph{local} action. But before we can do this, we need a precise statement about the definition of the functions we are dealing with.

\prop\label{prop:6.3a}
Let $f\in C_c(X)$ and $g\in C_c(G)$. Assume that the support of $g$ is contained in $V_f$. Define $F$ on $X\times G$ by
\begin{equation*}
F(x,p)=
\begin{cases}
f(x\cdot p)g(p) &\text{ if } (x,p)\in \Gamma\\
0 &\text{ otherwise.}
\end{cases}
\end{equation*}
Then $F$ is a continuous function with compact support on $X\times G$.
\eprop
\bew
i) The function $F$ is clearly continuous on $\Gamma$. 
\ssnl
ii) Denote by $C_f$ the support of $f$ in $X$ and by $C_g$ the support of $g$ in $G$. By assumption we have $C_g\subseteq V_f$. This means, by the definition  of $V_f$ (see Definition  \ref{defin:2.3}), that $(x,p\inv)\in \Gamma$ for $x\in C_f$ and $p\in C_g$. Then we can define $C$ in $X\times G$ as the set of points $(x\cdot p\inv,p)$ where $x\in C_f$ and $p\in C_g$. It is a compact subset of $\Gamma$.
\ssnl
iii) Now assume that $(x,p)\in \Gamma$ but $(x,p)\notin C$. If $x\cdot p$ is not in the support of $f$, we have $f(x\cdot  p)g(p)=0$. Suppose next that $(x\cdot p)\in C_f$. If $p\in C_g$ then, because $(x,p)=((x\cdot p)\cdot p\inv,p)$, we have $(x,p)\in C$. This contradicts the assumption. Therefore $p\notin C_g$. It follows that $g(p)=0$ and so also then $f(x\cdot p)g(p)=0$. 
\ssnl
iv) We find that  $F$ is $0$ on the complement of $C$. Because $C$ is a compact subset of $\Gamma$ and $F$ is continuous on $\Gamma$ we have $F$ continuous on all of $X\times G$. The support is contained in $C$ and so $F$ has compact support.
\ebew

Compare this result and its proof with Proposition \ref{prop:2.2}. 
Remark that $$F(x,p)=(p\cdot f)(x)g(p)$$ for all $x$ when $p\in V_f$. This follows from the fact that $(p\cdot f)(x)=f(x\cdot p)$ if $(x,p)\in \Gamma$ and $(p\cdot f)(x)=0$ if $(x,p)\notin\Gamma$, again see Proposition \ref{prop:2.2}.
\snl
We now generalize the result obtained of Proposition \ref{prop:6.1} to a local action.

\prop\label{prop:6.4a}
Assume that $f$ is a function in $C_c(X)$ and $g$ a function in $C_c(G)$ with support contained in $V_f$. Consider the function $F$ defined in the previous proposition. If $f$ is a polynomial function there are functions $g_i$ in $C_c(G)$ and  functions $f_i$ in $C_c(X)$  so that
\begin{equation}
F(x,p)=\sum_if_i(x)g_i(p)\label{eqn:6.2}
\end{equation}
for all all $x\in X$ and $p\in G$.  The functions $f_i$ can be chosen to be polynomial. Conversely, if we have a function $f\in C_c(X)$ with the property that there is a non-zero function $g\in C_c(G)$ with support contained in $V_f$ and functions $g_i$ in $C_c(G)$ and  functions $f_i$ in $C_c(X)$ such that Equation (\ref{eqn:6.2} holds, then $f$ must be polynomial. 
\eprop

\bew
i) Let $f$ be a polynomial function in $C_c(X)$. Take a function $g\in C_c(G)$ and assume that its support $C_g$ is contained in $V_f$. Now we choose a neighborhood $V$ of the identity in $G$ so that $V\subseteq V_f$. By taking it small enough, we can assume that also $C_gV\subseteq V_f$.  This is possible because $C_g$ is a compact set, contained in the open set $V_f$.  By taking $V$ still smaller if necessary, we can assume the existence of functions $f_i$ and $\varphi_i$  as in Proposition  \ref{prop:3.8} and we have $v\cdot f=\sum_i \varphi_i(v)f_i$ for all $v\in V$.
\ssnl
ii) Now we can proceed as in the proof of the global case. First choose elements $q_1, q_2,\dots,q_n$ in $C_g$ satisfying 
\begin{equation*}
C_g\subseteq q_1 V \cup q_2 V \cup \dots \cup q_n V.
\end{equation*} and a partition of unity $\lambda_1,\lambda_2,\dots\lambda_n$ of functions in $C_c(G)$, relative to this covering of $C_g$. Then, if $q\in V_f$  we can write
\begin{align*}
f(x\cdot q)g(q)
&=\sum_k \lambda_k(q)f(x\cdot q)g(q)\\
&=\sum_k \lambda_k(q) f(x\cdot q_k \cdot (q_k\inv q))g(q)\\
&=\sum_{i,k}  \lambda_k(q) f_i(x\cdot q_k)\varphi_i(q_k\inv q)g(q).
\end{align*}
However, because the support of the functions $\lambda_k$ also are contained in $V_f$, the equality will still be correct for the function $F$.
We have to take into account the following considerations. First we have that $(q_k\inv q)\cdot f$ is well-defined because $q_k\inv q\in V$ and $V\subseteq V_f$. Further we claim that $q_k\in V_{f_i}$ so that we can safely write
$q_k\cdot ((q_k\inv q)\cdot f)=q\cdot f$. To prove the claim we use the formula for $f_k$ as found in item iii) in the proof of Proposition \ref{prop:3.8}. We had
\begin{equation*}
f_i=\sum_j c_{ji}(v_{ji}\cdot f)
\end{equation*}
with $v_{ji}\in V$. But as $C_gV\subseteq V_f$ we know that $q_kv_{ji}\in V_f$ for all $j$ and therefore $q_k\in V_{f_i}$ for all $i,k$.
\ssnl
iii) Write $\psi_{ik}$ for the function $q\mapsto \lambda_k(q)\varphi_i(q_k\inv q) g(q)$. These are functions in $C_c(G)$ satisfying 
\begin{equation*}
F(x,q)=\sum_{ik} (q_k\cdot f_i)(x)\psi_{ik}(q)
\end{equation*}
as soon as $q\in V_f$. When $q\notin V_f$, not only $F(x,q)=0$, but also $\psi_{ik}(q)=0$ because all the functions $\lambda_k$ have support in $V_f$. Therefore the equation is still true in that case.
\ssnl
iv) To prove the converse result, we just can proceed as in the proof of item ii) of Proposition \ref{prop:6.1}.

\ebew
Also the  following corollary provides a characterization of polynomial functions.

\gev\label{gev:6.5} 
Let $f$ be a polynomial function on $X$. For any open subset $V$ of $G$ with compact closure contained in $V_f$, there are functions $g_i$ in $C_c(G)$ and polynomial functions $f_i$ in $C_c(X)$ so that $v\cdot f=\sum_i g_i(v)f_i$ for all $v\in V$.
\egev

\bew
Because the closure of $V$ is a compact subset of the open set $V_f$ we can find a function $g\in C_c(G)$ with compact support in $V_f$ and so that $g(v)=1$ for $v\in V$. If we  apply  Proposition \ref{prop:6.3a} we have functions $f_i\in C_c(X)$ and $g_i$ in $C_c(G)$ so that $g(q)(q\cdot f)=\sum_i g_i(q)f_i$ for all $q\in V_f$. If $v\in V$ then $g(v)=1$ and  we get $v\cdot f=\sum_i g_i(v)f_i$. This proves the result.
\ebew

Observe that in turn, the corollary implies the result of Proposition \ref{prop:6.4a}. Indeed, suppose that $f$ is a polynomial function on $X$ and that $g$ is a function in $C_c(G)$ with support in $V_f$.  Let $V$ be an open set containing the support of $g$ and such that its closure is a compact subset contained in $V_f$. Choose functions $g_i$ and $f_i$ as in the corollary satisfying $q\cdot f=\sum_i g_i(q)f_i$ for all $q\in V$. Now multiply this equation with $g(q)$, then we obtain
\begin{equation*}
g(q)(q\cdot f)=\sum_i g(q)g_i(v)f_i.
\end{equation*}
This holds in the first place for $q\in V$. But as the support of $g$ is contained in $V$, it holds for all $q\in V_f$.
\nl
\bf Polynomial functions on $\Gamma$ \rm
\nl
Given a local right action of $G$ on $X$, we consider the Cartesian product  $X\times G$. We have several actions of $G$ on this space. First we have a local right action given by the original action on the factor $X$ in this Cartesian product. Secondly we have left and right multiplication on the second factor. We have introduced these actions already in Section \ref{s:loc-act}, see Propositions \ref{prop:1.14} and \ref{prop:1.15}. We now consider polynomial functions for these actions.
\ssnl
Let $f\in C_c(X)$ and  $g\in C_c(G)$ with support in $V_f$. Consider the function $F$ on $X\times G$ defined by $F(x,p)=f(x\cdot p)g(p)$ when $(x,p)\in \Gamma$ and $F(x,p)=0$ otherwise (as in  Proposition \ref{prop:6.3a}).

\prop\label{prop:6.6}
If $f$ is polynomial, then $F$ is polynomial for the action of $G$ on the first factor of $X\times G$. If also $g$ is polynomial, then $F$ will be also polynomial for left and right multiplication by $G$ on the second factor in $X\times G$.
\eprop
\bew
i) We can write 
\begin{equation*}
F(x,p)=\sum_i f_i(x)g_i(p)
\end{equation*}
as we have seen in Proposition \ref{prop:6.4a}. The functions $f_i$ can be chosen to be polynomial. This implies that $F$ is polynomial for the action of $G$ on the factor $X$ in $X\times G$. 
\ssnl
ii) To prove that the function $F$ is polynomial in $G$ when this is the case for $g$, we proceed as follows.
\ssnl
We know that the support  $C_g$ of $g$ is contained in $V_f$. Then we can find a neighborhood $V$ of the identity in $G$, with compact closure $\overline V$ so that $C_g{\overline V}\inv\subseteq V_f$. We choose a function $\lambda\in C_c(G)$, with support in $V_f$ and  so that $\lambda(p)=1$ for $p\in C_g{\overline V}\inv$. Then we have $g(pv)=\lambda(p)g(pv)$ for all $p\in G$ and all $v\in V$.
\ssnl
iii) We now assume that $V$ is small enough so that it is  contained in $V_f$ and that we can write
$(v\cdot f)=\varphi_i(v)f_i$ for $v\in V$. Then we define functions $F_i$ on $X\times G$ by
\begin{equation*}
F_i(x,p)=
\begin{cases} 
f_i(x\cdot p)\lambda(p) &\text{ if } (x,p)\in \Gamma \\ 0 & \text{ otherwise.}
\end{cases}
\end{equation*} 
We know from Proposition \ref{prop:6.3a} that these functions are continuous with compact support on $X\times G$.
We claim that 
\begin{equation}
(v\cdot F)(x,p)=\sum_i \varphi_i(v) F_i(x,p) g(pv).\label{eqn:6.3}
\end{equation}
\ssnl
iv) To prove the claim, take an element $x\in X$ and consider elements $p\in G$ and $v\in V$. 
First suppose that $(x,pv)\in \Gamma$ and $(x,p)\in\Gamma$. Then 
\begin{align*}
(v\cdot F)(x,p)
&=F(x,pv)\\
&=f(x\cdot pv)g(pv)\\
&=f(x\cdot p\cdot v)g(pv)\\
&=(v\cdot f)(x\cdot p)g(pv)\\
&=\textstyle\sum_i \varphi_i(v)f_i(x\cdot p)\lambda(p)g(pv)\\
&=\textstyle\sum_i \varphi_i(v)F_i(x,p)g(pv).
\end{align*}
Next suppose that $(x,p)\notin\Gamma$. Then $F_i(x,p)=0$ for all $x$ and we have to argue that also $(v\cdot F)(x,p)=0$. Suppose however that $(v\cdot F)(x,p)\neq 0$. Then we must have $(x,pv)\in \Gamma$ and $F(x,pv)\neq 0$. We know that the support of $F$ is a compact subset of $\Gamma$. It must contain $(x,pv)$. By choosing $V$ small enough it follows that still $(x,p)\in \Gamma$. This gives a contradiction. Hence also $(v\cdot F)(x)=0$ and again  (\ref{eqn:6.3}) holds.
\ssnl
v) If now $g$ is polynomial and if we write $g(pv)=\sum_j \psi_j(v)g_j(p)$ we find in the end that 
\begin{equation*}
(v\cdot F)(x,p)=\sum_{i,j} \varphi_i(v)\psi_j(v) F_i(x,p) g_j(p)
\end{equation*}
for all $p$ and $v$ small enough. This proves that the function $p\mapsto F(x,p)$ is polynomial for all $x$. 
\ssnl
vi)
Finally, write again $F(x,p)=\sum_i f_i(x)g_i(p)$ and assume that the functions $f_i$ are linearly independent. Because $p\mapsto F(x,p)$ is polynomial for all $x$, we must have that $g_i$ is polynomial on $G$ for all $i$. This implies that the function $F$ itself is polynomial for the action of $G$ on the second factor in $X\times G$.
\ebew

\keepcomment{\rood Check again. What about the choice of $V$? Dependencies? See tmp file of 25 September 2020. Deze vind ik niet meer terug!
\nl}{}
In \cite{La-VD3} we illustrate this property with an example. \keepcomment {Complete with a concrete reference.}{}
\ssnl

This is a good place to refer again to Remark \ref{opm:3.10a}. A function on a Cartesian product can be polynomial in each factor, but need not be polynomial for the combined action.
\ssnl
At the end of the proof of the previous proposition, we showed that the function $F$ is polynomial for the action of $G$, from the fact that $p\mapsto F(x,p)$ is polynomial for each element $x$. This is possible because we can write $F(x,p)=\sum_i f_i(x)g_i(p)$.
\ssnl
This could be an indication to find a counter example for this problem.
\nl
\bf Polynomial functions and subgroups \rm
\nl
Assume that we have a local right action of a group $G$ on a space $X$ and that $H$ is a closed subgroup of $G$. Consider the restriction of the action of $G$ on $X$ to $H$ as in Definition \ref{defin:1.10}. 
\ssnl
\oldcomment{We need to change some conventions for notations. Sometimes we use $H$ for a closed subgroup, and sometimes for an open subgroup! This has been done - \rood but check here again!}{}
If $f$ is a polynomial function in $C_c(X)$ for the action of $G$, it is trivially also polynomial for the action of $H$, see  Proposition \ref{prop:3.3}. Then we can pose the following question.

\probl\label{probl:6.8a} 
Given a polynomial function $f$ in $C_c(X)$ for the action of the subgroup $H$. What  are possible extra conditions to ensure that the function $f$ is still polynomial for the original action of $G$\,? 
\eprobl
We will consider examples later, see Propositions \ref{prop:6.15a} and \ref{prop:6.18}).
\ssnl
Let us first recall what we know when $f$ is polynomial for the action of $G$. We have a compact subset $C$ of $X$, containing the support of $f$, an open subgroup $G_0$ of $G$, such that $C\cdot G_0$ is well-defined and equal to $C$, and a finite-dimensional subspace $F$ of the space  of functions in $C_c(X)$ with support in $C$, containing $f$, and invariant under the action of $G_0$. See Theorem \ref{stel:3.15}.  If the action is proper, the group $G_0$ is an open and compact subgroup of $G$. See Proposition  \ref{prop:3.18}. 
\ssnl
Further we now only assume that $f$ is a polynomial function on $X$ for the action of $H$. If we apply the above for this action, we find the following. \oldcomment{Formulate this and refer to it here! TO DO, see above.}

\prop\label{prop:6.9a} 
There is a compact subset $C$ of $X$, containing the support of $f$, and an open subgroup $H_0$ of $H$ such that $C\cdot H_0=C$. If we denote by $\mathcal F_C$ the space of functions $f\in C_c(X)$ with support in $C$, we have a representation $\pi$ of $H_0$ on $\mathcal F_C$ given by $\pi(h)f)(x)=f(x\cdot h)$. Finally $f$ belongs to a finite-dimensional $\pi$-invariant subspace of $\mathcal F_C$.
\eprop
Remark that we should rather write $h\cdot f$ for $\pi(h)f$, but as we know that the support of $f$ is contained in $C$ and that $C$ is invariant under the action of $H_0$, we will have $(\pi(h)f)(x)=f(x\cdot h)$ whenever this is non-zero.
\ssnl
It is clear that nothing can be said further without extra conditions. Think of the possibility that the action of $H$ is trivial. 
\ssnl
On the other hand, with the following condition, we can proceed in the direction of the desired property.

\voorw\label{voorw:6.10a} 
In what follows, we assume that the action of $H$ on $X$ is locally homeomorphic as in Definition \ref{defin:1.6}. So for each $x\in X$ we have a neighborhood $V_x$ of the identity in $H$ so that $x\cdot V_x$ is a well-defined open subset of $X$ and so that the map $h\mapsto x\cdot h$ is a homeomorphism from $V_x$ to $x\cdot V_x$. 
\evoorw

Because $H_0$ is an open subgroup and $h\mapsto x\cdot h$ is a homeomorphism from $V_x$ to its image, we can replace $V_x$ with the intersection $H_0\cap V_x$ and still have the same property. So in what follows we assume that $V_x\subseteq H_0$.  
\ssnl
We now continue under the assumption that Condition \ref{voorw:6.10a} holds and use the notation of Proposition \ref{prop:6.9a} above. Then we clearly get the following.

\prop\label{prop:6.11a} 
The set $C$ is a compact open subset of $X$.
\eprop

\bew
We know that $C$ is a compact subset of $X$. Moreover, for all $x\in X$ we have the neighborhood $x\cdot V_x$ of $x$ in $X$. Because we assume that $V_x\subseteq H_0$, we have that  $x\cdot V_x\in C$. Therefore $C$ is also open.
\ebew

\prop\label{prop:6.12a} 
Define $G_1$ as the subset of points $p$ in $G$ satisfying $(x,p)\in\Gamma$ and $x\cdot p\in C$ for all $x\in C$. Then $G_1$ is an open subset of $G$ and $pq\in G_1$ if $p,q\in G_1$. 
\eprop
\bew
i) Because $C$ is open, for all $x\in C$ we have a neighborhood $U_x$ of the identity in $G$ so that $x\cdot U_x$ is well-defined and contained in $C$. By a standard argument, because $C$ is compact, we find a single neighborhood $U$ of the identity in $G$ so that $x\cdot U$ is defined and contained in $C$ for all $x\in C$, see e.g.\ the argument in the proof of Proposition \ref{prop:2.3}. \oldcomment{I am not sure if that is a good internal reference.}{} 
Therefore $U\subseteq G_1$.
\ssnl
ii) Now assume that $p,q\in G_1$. Then
\begin{equation*}
C\cdot (pq)=(C\cdot p)\cdot q\subseteq C\cdot q\subseteq C
\end{equation*}
so that also $pq\in G_1$. 
\ssnl
iii) Given $p\in G_1$ we have a neighborhood $pU$ of $p$ in $G$, contained in $G_1$. So $G_1$ is open in $G$.
\ebew
We can not conclude at this stage that $G_1$ is actually a group, but it certainly contains an open subgroup $G_0$ if we take the neighborhood $U$ symmetric and the group of elements generated by $U$. On the other hand, $H_0$ is contained in $G_1$. And if we replace $H_0$ by $H\cap G_0$, we still get an open subgroup of $H$. This will eventually lead to the following refinement of 
Theorem \ref{stel:3.15} of Section \ref{s:pol}.

\prop\label{prop:6.13a} 
We have an open and compact subset $C$ of $X$ and an open subgroup $G_0$ of $G$ so that $C\cdot G_0=C$. This yields a representation of $G_0$ on the space $\mathcal F_C$ of functions in $C_c(X)$ with support in $C$. If we let $H_0=G_0\cap H$, we have an open subgroup $H_0$ of $H$, contained in $G_0$. The function $f$ belongs to a finite-dimensional $H_0$ invariant subspace of $\mathcal F_C$.
\eprop

We will show later that a finite-dimensional representation of $H_0$ can occur at most  finitely many times in this representation of $G_0$, see Proposition \ref{prop:6.16a} below. The following property is related to that.

\prop\label{prop:6.14a} 
As before, we have a compact subset $C$ of $X$ and we use $\mathcal F_C$ for the space of functions in $C_c(X)$ with support contained in $C$. 
Suppose that we have finitely many functions $\varphi_i$ on $H_0$. Consider the subspace $F$ of functions $f$ in $\mathcal F_C$ such that there exist functions $g_i $ in $\mathcal F_C$ with the property that $f(x\cdot h)=\sum_i \varphi_i(h)g_i(x)$ for all $x\in C$ and all $h\in H_0$. Then $F$ is finite-dimensional.
\eprop
\bew
For each element $x\in C$ we choose a neighborhood $V_x$ of the identity of $H$. We can assume that $V_x\subseteq H_0$. We can cover $C$ by finitely many sets $x_kV_k$ where we use $V_k$ for the neighborhood $V_{x_k}$. We also denote by $p_k$ the continuous function from $x_k\cdot V_k$ determined by $x=x_k\cdot p_k(x)$ for $x\in x_k\cdot V_k$. Consider a partition of unity $(\lambda_k)$ relative to this covering.
\ssnl
Now let $g$ be in $F$. Then we have
\begin{align*}
g(x)
&=\sum_k \lambda_k(x) g(x)\\
&=\sum_k \lambda_k(x) g(x_kp_k(x)) \\
&=\sum_{j,k} \lambda_k(x) g_i(x_k)\varphi_i(p_k(x)).
\end{align*}
We see that $g$ belongs to the space spanned by the functions $x\mapsto \lambda_k(x)\varphi_i(p_k(x))$. This proves the result.
\ebew

This allows us to obtain a solution for Problem \ref{probl:6.8a}.

\prop\label{prop:6.15a} 
Let $G$ be a locally compact group with a compact open subgroup and  let $H$ be a closed subgroup of $G$.  Assume that we have a local right action of $G$ on $X$ and that the restriction of the action to $H$ is locallly homemorphic. Then any function $f\in C_c(X)$ that is polynomial for the action of $H$ will still be polynomial for the action of the larger group $G$.
\eprop

\bew
Let $f$ be a polynomial function in $C_c(X)$ for the action of $H$. We apply Proposition \ref{prop:6.13a}.  We find a compact open subset $C$ of $X$ that contains the support of $f$, an open subgroup $G_0$ of $G$ so that $C\cdot G_0=C$ such that $f$ belongs to a finite-dimensional subspace $F$ of $\mathcal F_C$ that is invariant under the action of the open subgroup $H_0=G_0\cap H$. As we now assume that $G$ has a compact open subgroup, we can take the intersection of $G_0$ with this compact open subgroup and assume in the end that $G_0$ itself is a compact open subgroup of $G$.
\ssnl
We have functions $g_i\in\mathcal F_C$ and $\varphi_i\in C_c(G)$ satisfying $(h\cdot f)=\sum_i \varphi_i(h)g_i$ for all $h\in H_0$. Take any irreducible representation $\gamma$ of $G_0$ and consider the projection map $P_\gamma$ on the space $\mathcal F_C$. Because $H_0$ is a subgroup of $G_0$, we have $\pi(h)P_\gamma=P_\gamma\pi(h)$ so that still $h\cdot (P_\gamma f)= \sum_i \varphi_i P_\gamma g_i$. By Proposition \ref{prop:6.14a} all the functions $P_\gamma f$ belong to a finite-dimensional subspace for all $\gamma$. This means that only finitely many of them can be non-zero. This implies that $f$ is polynomial for the action of $G_0$.
\ebew

We will use this property later in this paper, see Proposition \ref{prop:6.18}. \keepcomment {Also in \cite{La-VD3} - References?
\ssnl}

We now include some more information to get a better understanding of this result.
 
\prop\label{prop:6.16a} 
As in the previous proposition, let $C$ be a compact  subset of $X$ and assume that $H_0$ is a compact group acting on $X$ and satisfying $C\cdot H_0=C$. Assume that the action is locally homeomorphic. Consider the associated representation $\pi$ of $H_0$ on the space $\mathcal F$ of functions in $C_c(X)$ with support in $C$. Then   
every irreducible representation $\gamma$ of $H$ can occur at most finitely many times in $\pi$.
\eprop
\bew
Because the action of $H_0$ is assumed to be locally homeomorphic, we can cover $C$ by finitely many, say $n$ open sets $x_kH_0$.  Define a map $T$ from $\mathcal F$ to the direct sum of $n$ copies of $C(H_0)$ by $(Tf)_k(h)=f(x_kh)$. The map is injective because if $f(x_kh)=0$ for all $k$ and all $h$, then $f(x)=0$ for all $x\in C$. The map $T$ is also an intertwiner between the representation $\pi$ on $\mathcal F$ and the direct sum of $n$ copies of the right regular representation of $H_0$ on $C(H_0)$. Now because each irreducible representation of $H_0$ can occur only finitely many times in the regular representatiion, the result follows.
\ebew


We have some applications of this result in \cite{La-VD3}
\nl
\bf Polynomial functions for commuting actions \rm
\nl
In Proposition \ref{prop:6.18} below we obtain an application of the result in Proposition \ref{prop:6.15a} for commuting actions. First we prove some general other result about polynomial functions for such commuting actions.
\ssnl
Suppose that we have a local left  action of a group $H$  and a local right action of a group $K$ on $X$. Assume that they commute as in Definition \ref{defin:1.2}. Consider a function $f\in C_c(X)$ that is polynomial for the associated right local action of $H^\text{op}\times K$. Then it is polynomial for the left action of $H$ and it is polynomial for the right action of $K$. This is so because these actions are restrictions to a closed subgroup, see Definition \ref{defin:1.10}. 
\ssnl
The converse is also true. It is not completely obvious as we see from the proof below.

\prop\label{prop:6.17a}
Consider a local left action of $H$ and a local right action of $K$ on $X$. Assume that they commute as in Definition \ref{defin:1.2}. Let $f$ be in $C_c(X)$ and assume that it is polynomial both for the action of $H$ and for the action of $K$. Then it is polynomial for the associated local right  action of $H^\text{op}\times K$.
\eprop
\bew
The associated action of  $H^\text{op}\times K$ is considered in  a comment following Proposition \ref{prop:1.3}.
\ssnl
i) First we choose neighborhoods $U$ and $V$ of the identity in $H$ and $K$ respectively so that $v\cdot f\cdot u$ is well-defined for $u\in U$ and $v\in V$ as 
\begin{equation*}
(v\cdot f\cdot u)(x)=f(u\cdot x\cdot v)
\end{equation*}
when $u\cdot x\cdot v$ is defined (cf.\ Lemma \ref{lem:2.9a}). 
\ssnl
ii) Next we use that $f$ is polynomial for the right action of $K$. Then we have  functions $\varphi_i$ and $f_i$  so that we can write
\begin{equation*}
v\cdot f=\sum_i \varphi_i (v) f_i
\end{equation*}
for all $v\in V$. As in the proof of Proposition \ref{prop:3.8}, we can use Lemma \ref{lem:3.5} to  find numbers $c_{jk}$ and elements $v_{jk}\in V$ so that, for all $k$ we have 
\begin{equation*}
f_k=\sum_j c_{jk} (v_{jk}\cdot f).
\end{equation*}
\ssnl
iii) Now we use that $f$ is polynomial for the left action of $H$. Then we have functions $\psi_j$ and $f'_j$ so that we can write, for all $u\in U$,
\begin{equation*}
f\cdot u=\sum_j \psi_j(u)f'_j.
\end{equation*}
iv) Then we have 
\begin{align*}
v\cdot f\cdot u
&=\sum_i \varphi_i(v)(f_i\cdot u) \\
&=\sum_{i,k}  c_{ik}\, \varphi_i(v)(v_{ik}\cdot f\cdot u) \\
&=\sum_{i,j,k}  c_{ik}\, \varphi_i(v)\psi_j(u) (v_{ik}\cdot f'_j).
\end{align*}
First look at the second equality. We use that $v_{ik}\cdot f\cdot u$ is well-defined. This implies that $f_i\cdot u$ is well-defined. And we have the first equality. For the last equality we need that $v_{ik}\cdot f'_j$ is well defined. This is the same argument as above for showing that $f_i\cdot u$ is well-defined.
\ssnl
It follows that there is a finite-dimensional subspace of $C_c(X)$ that contains all the functions $v\cdot f\cdot u$. This shows that $f$ is polynomial for the joint action.
\ebew 

We would like to make the following remark in connection with the above result.

\opm\label{opm:6.17}
i) Suppose we have locally compact spaces $X$ and $Y$ with a local right action of a locally compact group $H$ on $X$ and a local right action of a locally compact group $G$ on $Y$. They yield commuting actions on $X\times Y$. Now let $f$ be a function in $C_c(X\times Y)$ and assume that it is polynomial for the action of $H$, as well for the action of $K$ on $X\times Y$. Then one can wonder if it will be polynomial for the combined action of $H\times K$. The problem is of a similar nature as the one considered earlier in Remark \ref{opm:3.10a}. 
\ssnl
ii) As a special case, we can consider two locally compact groups $H$ and $K$ and a function $f$ in $C_c(H\times K)$. Then the question is, if given that $f(\,\cdot\, ,k)$ is polynomial on $H$ for each $k$ and that $f(h,\,\cdot\,)$ is polynomial on $K$ for each $h$, does it follow that $f$ is a polynomial function on $H\times K$.
\ssnl
iii) The above proposition is of a similar nature, but still fundamentally different. \keepcomment{\rood Explain! - Still to be done.}{}
\eopm

Now we come to the application of Proposition \ref{prop:6.15a} for commuting actions.
  
\prop\label{prop:6.18}
Assume that we have a right local action of $G$ that is locally homeomorphic and a polynomial function for this action. Then $f$ will be polynomial for any left local action of a group $H$ that commutes with the given right local action.
\eprop
\bew
i) Because $f$ is polynomial, we have a neighborhood of the identity $V$ such that $V$ is contained in $V_f$ and we can write
\begin{equation*}
v\cdot f=\sum_i\varphi_i(v)f_i.
\end{equation*}
\ssnl
ii) Now use $(h,x)\mapsto h\cdot x$ for the left local action of $H$. As we have done for right actions, for a function $g$ in $C_c(X)$ and an element $u\in H$ we define the function $g\cdot u$ on $X$ by
\begin{equation*}
(g\cdot u)(x)=
\begin{cases}
g(u\cdot x) & \text{ if } u\cdot x \text{ is defined}\\
0 & \text{ otherwise.}
\end{cases}
\end{equation*}
For $u$ in a small enough neighborhood $U$ of the identity in $H$  depending on $g$, this function will be in $C_c(X)$. The argument is as the one given for the function $v\cdot f$ in Section \ref{s:ind-fun}. Then we can safely write
$(g\cdot u)(x)$ as $g(u\cdot x)$. 
\ssnl
iii) By taking $U$ sufficiently small we have continuous functions $f_i\cdot u$ for all $i$ where the functions $f_i$ are as in item i) of the proof above. Then we will have 
\begin{equation*}
v\cdot f\cdot u=\sum_i \varphi_i(v)f_i\cdot u 
\end{equation*}
for $u\in U$. 
\ssnl
iv) If we take for  $U$ a neighborhood with a compact closure, the supports of the functions $f\cdot u$ will all lie in a single compact subset of $X$ and so we can apply Lemma \ref{lem:2.9a}. We find that the functions $f\cdot u$ all belong to a finite-dimensional subspace of $C_c(X)$. Therefore $f$ is polynomial for the left action of $H$ on $X$
\ebew

We have some similar result in Proposition \ref{prop:5.12}. In fact this result is a special case of the above proposition because right multiplication by $G$ on $G$ is locally homeomorphic.
\ssnl

\section{\hspace{-17pt}. Conclusions and further research}\label{s:concl} 

In this note, we have treated polynomial functions for a partial action of a locally compact group $G$ on a locally compact space. We have shown that for proper actions, non-trivial polynomial functions only can exist when the group $G$ has a compact open subgroup  and there are compact subsets invariant under the action of this subgroup. 
\ssnl
In an earlier paper \cite{La-VD1} we have defined and studied polynomial functions on a locally compact group. This is a special case of the situation in this note. We recover some of the results obtained in \cite{La-VD1}. However, the treatment is different in the sense that here we define polynomial functions from the very beginning without a reference to a compact open subgroup as it was done in \cite{La-VD1}. 
\ssnl
This has some advantages. In particular, it allows us to generalize the notion and define polynomial functions on a locally compact space  relative to only a local action of a group. Because the action is not a global action, some care is needed. But as it turns out, polynomial functions in this setting behave very much like polynomial functions for a global action, in particular for the action of a group by multiplication on itself.
\ssnl
This generalization of polynomial functions to local actions is important for  its application to the theory of bicrossproducts for locally compact groups with a compact open subgroup as  treated in \cite{La-VD3}.
\ssnl
By assumption, polynomial functions on a locally compact group $G$, or more generally on a locally compact space $X$ with a local action of a locally compact group, are required to have compact support. If we drop this requirement, and simply require to have continuous functions, the situation is completely different. But it still makes sense to investigate this notion.  
\ssnl
We also feel that some more research should be done. In this note, we had the applications in mind that we need for the special cases we use in \cite{La-VD3}. That is the reason why we have included such examples in Section \ref{s:misc}. But other examples to illustrate our results would be most welcome. We already mentioned a problem in Remark \ref{opm:3.10a} and refered to this in Section \ref{s:misc}, after the proof of Proposition \ref{prop:6.6}, with an indication for how to find an example to solve this problem.
\ssnl
Another issue is the problem to determine which local actions are obtained by restriction of a global action to a subspace as in Proposition \ref{prop:1.8}. Most of the examples we have treated in this note are of this form, but some others are not, at least not obviously. See the remark following the proof of Proposition \ref{prop:1.15}.
We believe that the problem is open and interesting. Remark however that this has nothing to do with the theory of polynomial functions. 
\ssnl
Finally, we refer to \cite{La-VD3} where local actions and polynomial functions are used in the development of algebraic quantum hypergroups arising from locally compact groups with a compact open subgroup.
\nl



\end{document}